\documentclass[a4paper,12pt,leqno]{article}
\usepackage[latin1]{inputenc}
\usepackage{NVWMacro}
\usepackage{amsmath}
\usepackage{amsfonts}
\usepackage{amssymb}
\usepackage[all]{xy}
\usepackage{mathrsfs}
\usepackage{stmaryrd}
\usepackage{bbm}
\usepackage{oldgerm}
\usepackage{pst-all,multido,ifthen}

\newcommand{\Cond}{({\rm C})}

\begin{document}

\title{Purity of level $m$ stratifications}

\author{Marc--Hubert Nicole\footnote{nicole@math.jussieu.fr} , Adrian  
Vasiu\footnote{adrian@math.binghamton.edu} , Torsten  
Wedhorn\footnote{wedhorn@math.uni-paderborn.de} }

\maketitle


\noindent{\scshape Abstract.\ } Let $k$ be a field of
characteristic $p>0$. Let $D_m$ be a $\BT_m$  over $k$ (i.e., an
$m$-truncated Barsotti--Tate group over $k$). Let $S$ be a\break
$k$-scheme and let $X$ be a $\BT_m$ over $S$. Let $S_{D_m}(X)$ be
the subscheme of $S$ which describes the locus where $X$ is
locally for the fppf topology isomorphic to $D_m$. If $p\ge 5$, we
show that $S_{D_m}(X)$ is pure in $S$ i.e.,~the immersion
$S_{D_m}(X) \hookrightarrow S$ is affine. For $p\in\{2,3\}$, we
prove purity if $D_m$ satisfies a certain technical property depending only
on its $p$-torsion $D_m[p]$. For $p\ge 5$, we apply the developed
techniques to show that all level $m$ stratifications associated
to Shimura varieties of Hodge type are pure.

\

\noindent{\scshape R\'esum\'e.\ } Soit $k$ un corps de
caract\'eristique $p>0$. Soit $D_m$ un $\BT_m$ sur $k$ (i.e., un
groupe de Barsotti--Tate tronqu\'e en \'echelon $m$ sur $k$). Soient
$S$ un $k$-sch\'ema et $X$ un $\BT_m$ sur $S$. Soit $S_{D_m}(X)$
le sous-sch\'ema de $S$ correspondant au lieu o\`u $X$ est
isomorphe \`a $D_m$ localement pour la topologie fppf. Si $p\ge
5$, nous montrons que $S_{D_m}(X)$ est pur dans $S$ i.e.,
l'immersion $S_{D_m}(X) \hookrightarrow S$ est affine. Pour
$p\in\{2,3\}$, nous prouvons la puret\'e pour $D_m$ satisfaisant une
certaine propri\'et\'e technique d\'ependant uniquement de la $p$-torsion
$D_m[p]$. Pour $p \ge 5$, nous utilisons les techniques
d\'evelopp\'ees pour montrer que toutes les stratifications par l'\'echelon associ\'ees aux vari\'et\'es de Shimura de type Hodge
sont pures.




\section{Introduction}

Let $p$ be a prime number. Let $k$ be a field of
characteristic $p$. Let $c$, $d$, and $m$ be positive integers. In
this paper, a $\BT_m$ is an $m$-truncated Barsotti--Tate group of
codimension $c$ and dimension $d$. Let $D_m$ be a fixed $\BT_m$
over $k$.

Let $S$ be an arbitrary $k$-scheme and let $X_m$ be a $\BT_m$ over
$S$. Let $S_{D_m}(X_m)$ be the (necessarily unique) locally closed subscheme of $S$ that satisfies the following property. A morphism $f\colon S' \to
S$ of $k$-schemes factors through $S_{D_m}(X_m)$ if and only if
$f^*(X_m)$ and $D_m \times_{\Spec k} S'$ are locally for the fppf
topology isomorphic as $\BT_m$'s over $S'$ (see
Subsection~\eqref{LevelmStrata} for the existence of
$S_{D_m}(X_m)$). If $D$ is either a $\BT_{m'}$ for some $m' \geq m$,
or a $p$-divisible group over $k$, we will also write $S_D(X_m)$
instead of $S_{D[p^m]}(X_m)$.

The following notion of purity (that has already been considered in  
\cite{Va1}, Section~2.1.1) will be central.

\begin{proclamation}{Definition}\label{DefPure}
A subscheme $T$ of a scheme $S$ is called \emph{pure in $S$} if the immersion $T\hookrightarrow S$ is affine.
\end{proclamation}

We remark that the purity of $T$ in a locally noetherian scheme
$S$ implies the following weaker variant of purity: If $Y$ is an
irreducible component of the Zariski closure $\Tbar$ of
$T$ in $S$, then the complement of $Y \cap T$ in $Y$ is either
empty or of pure codimension $1$. On the other hand, if $S$ is separated and $T$ is (globally) an affine scheme, then $T$ is pure in $S$.

Purity results for strata defined by $p$-divisible groups have a
long history. The earliest hints of purity are probably the
computations mentioned by Y.~Manin in \cite{Ma}, at the bottom of
p. 44. For Newton polygon strata, J.~de Jong and F.~Oort have
shown the above mentioned weaker version of purity in \cite{dJO}
and one of us has shown in \cite{Va1} that these strata are even
pure in the sense of Definition~\ref{DefPure}. For $p$-rank
strata, Th.~Zink proved in \cite{Zi} the weaker version of purity.
Moreover, T.~It\=o proved in \cite{Ito} the existence of
generalized Hasse--Witt invariants for PEL unitary Shimura
varieties of signature $(n-1,1)$ at primes $p$ where the unitary
group is split. This result implies in fact a stronger kind of
purity (see below).

The weak version of purity is an important tool to estimate and
compute the dimensions of strata in the locally noetherian case.
Purity itself is an important step towards determining whether a
(quasi-affine) stratum is in fact affine, or whether a cohomological sheaf is in fact zero. For instance, a genuine (cohomological) application of purity (and not of global affineness!) to some simple Shimura varieties can be found in 
\cite{Bo}, Proposition 6.2.

The goal of this paper is to show that $S_{D_m}(X_m)$ is pure for all  
schemes $S$ and all $BT_m$'s $X_m$ if $D_m$ satisfies a certain  
condition $\Cond$ introduced in Subsection~\eqref{COND}. Here we  
remark that condition $\Cond$ depends only on $D_m[p]$ and it can be  
checked easily. Condition $\Cond$ is satisfied if any one of the  
three conditions below holds (cf. Lemma~\ref{CompLem}~(c)
and~(d) and Example \ref{Ex1}):
\begin{romanlist}
\item We have $p\ge 5$.
\item We have $p = 3$ and $\min(c,d) \le 6$.
\item There exists an integer $a\ge 2$
such that we have a ring monomorphism $\Bbb
F_{p^a}\hookrightarrow \End(D_m[p])$ with the property that $\Bbb
F_{p^a}$ acts on the tangent space of $D_m[p]$ via scalar
endomorphisms.
\end{romanlist}

\medskip
\noindent For the remainder of the introduction, we will assume
that condition $\Cond$ holds for $D_m$. The main result of the paper  
is the following
theorem.

\begin{proclamation}{Theorem}\label{MainThm}
The locally closed subscheme $S_{D_m}(X_m)$ is pure in $S$.
\end{proclamation}

We obtain the following corollary.

\begin{proclamation}{Corollary}
Let $S$ be locally noetherian and let $Y$ be an irreducible
component of $\overline{S_{D_m}(X_m)}$. Then
the complement of $S_{D_m}(X_m) \cap Y$ in $Y$ is either empty or of
pure codimension $1$.
\end{proclamation}

Now let $D$ be a $p$-divisible group over $k$ such that $D[p^m]=D_m$.  
For every reduced $k$-scheme $S$ and every $p$-divisible group $X$  
over $S$ denote by $\grn_D(X)$ the (necessarily unique) reduced  
locally closed subscheme of $S$ such that for each field extension $K$  
of $k$ we have
\[
\grn_D(X)(K) = \set{s \in S(K)}{\text{$D$ and $s^*(X)$ have equal  
Newton polygons}}.
\]
Thus $\grn_D(X)$ is the Newton polygon stratum of $S$ defined by $X$  
that corresponds to the Newton polygon of $D$. The locally closed subscheme  
$\grn_D(X)$ is pure in $S$ by~\cite{Va1}, Theorem~1.6. Thus we get  
another purity result:

\begin{proclamation}{Corollary}
For each $m\in\Bbb N^*$, the locally closed subscheme $\grn_D(X) \cap S_D(X[p^m])$ is pure in $S$.
\end{proclamation}

Moreover, we can use the well known fact that there exists an integer $n_D \geq  1$ with the following property. If $C$ is a $p$-divisible group over  
an algebraic closure $\kbar$ of $k$ such that $C[p^{n_D}]$ is isomorphic to  
$D[p^{n_D}]_{\kbar}$, then $C$
is isomorphic to $D_{\kbar}$ (for instance, see \cite{Tr1}, Theorem~1 or
\cite{Va1}, Corollary~1.3). We assume that $n_D$ is chosen minimal.  
Then there exists a (necessarily unique) reduced locally closed  
subscheme $\ufr_D(X)$ of $S$ such that for every algebraically closed  
field extension $K$ of $k$ we have
\[
\ufr_D(X)(K) = \set{s \in S(K)}{D_K \cong s^*(X)}.
\]
Indeed, we have $\ufr_D(X) = S_D(X[p^{n_D}])_{\rm red}$. From
Theorem~\ref{MainThm}, we obtain the following purity result:

\begin{proclamation}{Corollary}\label{UltStrat}
The locally closed subscheme $\ufr_D(X)$ is pure in $S$.
\end{proclamation}

For special fibres of good integral models in unramified mixed
characteristic $(0,p)$ of Shimura varieties of Hodge type (or more
generally, for quasi Shimura $p$-varieties of Hodge type), there
exists a level $m$ stratification that parametrizes
$\BT_m$'s with additional structures (see
Subsection~\eqref{QuasiShimura}). The proof of
Theorem~\ref{MainThm} can be adapted to show that all level
$m$ stratifications are pure (see Theorem \ref{PurityShimura}),
provided they are either in characteristic $p\ge 5$ or are in
characteristic $p\in\{2,3\}$ and an additional condition holds.

In this introduction, we will only state the Siegel modular varieties variant of
Theorem \ref{MainThm} (see Example \ref{Ex5}). Let $N\ge 3$ be an  
integer prime to $p$. Let $\Ascr_{d,1,N}$ be the Mumford moduli scheme  
that parameterizes principally polarized abelian schemes over  
$\FF_p$-schemes of relative dimension $d$ and equipped with a  
symplectic similitude level $N$ structure (cf.~\cite{MFK}, Theorems  
7.9 and 7.10). Let $(\Uscr,\Lambda)$ be the principally  
quasi-polarized $p$-divisible group of the universal principally  
polarized abelian
scheme over $\Ascr_{d,1,N}$. If $k$ is algebraically closed and if  
$(D,\lambda)$ is a principally
quasi-polarized $p$-divisible group of height $2d$ over $k$, let
$\grs_{D,\lambda}(m)$ be the unique reduced locally closed subscheme  
of $\Ascr_{d,1,N,k}$ that satisfies the following identity of sets
\[
\grs_{D,\lambda}(m)(k) =  
\set{y\in\Ascr_{d,1,N}(k)}{y^*(\Uscr,\Lambda)[p^m] \cong  
(D,\lambda)[p^m]}.
\]
Then $\grs_{D,\lambda}(m)$ is regular and equidimensional (see \cite{Va2}, Corollary 4.3 and Example 4.5; Subsection~\eqref{Dev} below can be easily adapted to prove the existence and the smoothness of the $k$-scheme $\grs_{D,\lambda}(m)$). Moreover we have:

\begin{proclamation}{Theorem}\label{SiegelThm}
If either $p=3$ and $d\le 6$ or $p\ge 5$, then the locally closed subscheme  
$\grs_{D,\lambda}(m)$ is pure in $\Ascr_{d,1,N,k}$.
\end{proclamation}

We remark that for $m=1$, Theorem~\ref{SiegelThm} neither implies nor  
is implied by Oort's result (\cite{Oo2}, Theorem~1.2) which asserts  
(for all primes $p$) that the scheme $\grs_{D,\lambda}(1)$ is  
quasi-affine.

Finally, we investigate briefly the following stronger notion of purity.

\begin{proclamation}{Definition}\label{DefPrincPure}
Let $T \to S$ be a quasi-compact immersion and let $\bar T$ be the
scheme-theoretic closure of $T$ in $S$. Then $T$ is called
\emph{Zariski locally principally pure in $S$} if locally for the
Zariski topology of $\bar T$, there exists a function $f \in
\Gamma(\bar T, \Oscr_{\bar T})$ such that we have $T = \bar T_{f}$, where $\bar T_{f}$ is the largest open subscheme of $\bar T$ over  
which $f$ is invertible.
\end{proclamation}

We obtain variants of this notion by replacing the Zariski
topology by another Grothendieck topology $\Tscr$ of $S$. If
$\Tscr$ is coarser than the fpqc topology (e.g., the Zariski or
the \'etale topology), each $\Tscr$ locally principally pure
subscheme is pure (as affineness for morphisms is a local property for the fpqc
topology). Principal purity for $p$-rank strata corresponds to the
existence of generalized Hasse--Witt invariants. They have been
investigated by T. It\=o for certain unitary Shimura varieties (see  
\cite{Ito}) and by E. Z. Goren for Hilbert modular varieties (see  
\cite{Go}).

In Section~\ref{StrongPurity}, we will show that this stronger
notion of purity does not hold in general. In fact, we have:

\begin{proclamation}{Proposition}\label{NoPrincPure}
Let $c,d \geq 2$ and $s\in \{1,\ldots,c-1\}$. Then the strata of $p$-rank equal to $s$ associated to $BT_1$'s over $\Bbb F_p$-schemes of codimension $c$ and dimension $d$, are not \'etale locally principally pure in general.
\end{proclamation}

\noindent
We now give an overview of the structure of the paper. In
Section~\ref{devissage}, we define the level $m$ strata
$S_{D_m}(X_m)$ and we prove some basic properties of them. Then we make a
d\'evissage to the following situation.

\begin{remark}{Essential Situation}\label{EssSit}
Let $k$ be an algebraically closed field of characteristic $p>0$
and let $D_m$ be a $\BT_m$ over $k$ which
satisfies condition~$\Cond$. Let $D$ be a $p$-divisible group over $k$
such that $D[p^m]=D_m$. Let $S = \Ascr$ be a smooth $k$-scheme of
finite type which is equidimensional of dimension $cd$ and for
which the following two properties hold:

\medskip
\begin{definitionlist}
\item
There exists a $p$-divisible group $\Escr$ of codimension $c$ and  
dimension $d$ over
$\Ascr$ which is a versal deformation at each $k$-valued point of $\Ascr$.
\item
There exists a point $y_D\in\Ascr(k)$ such that $y_D^*(\Escr)$ is  
isomorphic to $D$.
\end{definitionlist}
\end{remark}

\noindent In this case we simply write $\grs_{D}(m)$ instead of
$\Ascr_D(\Escr[p^m])$. In Subsection~\eqref{Dev}, we will prove that  
$\grs_{D}(m)$ is smooth over $k$ (by \cite{Va2}, Theorem~1.2~(a)
and~(b) and Remark~3.1.2 we know already that the reduced scheme of  
$\grs_{D}(m)$ is a smooth equidimensional $k$-scheme, although this  
fact is not used in the proof below). Then we show that  
Theorem~\ref{MainThm} follows if $\grs_D(m)$ is pure in $\Ascr$.

We remark that for $m \geq n_D$ (where $n_D$ is the integer defined  
above before Corollary~\ref{UltStrat}) the fact that $\grs_D(m)$ is  
pure in $\Ascr$ is proved in \cite{Va1}, Theorem~5.3.1 (c). This result of \cite{Va1} and thus Corollary \ref{UltStrat} also, hold even if 
condition $\Cond$ does not hold for $D[p]$.

The proof of Theorem~\ref{MainThm} is presented in Section~\ref{MainProof}. There we show that purity follows from the affineness of a certain orbit $\Oscr_m$ of a group action
\[
\Bbb T_m\colon \Hscr_m \times \Dscr_m \to \Dscr_m,
\]
which was introduced in \cite{Va2}. The orbits of $\Bbb T_m$ parameterize isomorphism classes of $\BT_m$'s over perfect fields. In fact we show that $\Oscr_m$ is affine for all $m$ provided $\Oscr_1$ is affine.

The definition of the action $\Bbb T_m$ is recalled in Section~\ref{GroupActions}, and in Section~\ref{Combinatorics} the main properties of the action $\Bbb T_1$ we
need are presented. There the condition~$\Cond$ is introduced and used. Its key role is to imply
that a certain morphism between affine $k$-schemes is finite (see Theorem~\ref{FinThm}) which allows us in Subsection ~\eqref{ProofofMain} to use Chevalley's theorem to show that $\Oscr_1$ is affine.

In Section~\ref{StratHodge}, we prove purity of the level $m$
stratification for quasi Shimura $p$-varieties of Hodge type.
Finally, in Section~\ref{StrongPurity}, we show that generalized
Hasse--Witt invariants do not always exist for Shimura varieties.
Indeed, we construct a counterexample and we deduce
Proposition~\ref{NoPrincPure}.


\section{D\'evissage to the Essential Situation}\label{devissage}

\begin{segment}{Moduli spaces of truncated Barsotti--Tate groups}{ModTruncBT}
We recall that $c$, $d$, and $m$ are positive integers. Let
$\Bscr\Tscr_m = \Bscr\Tscr_m^{c,d}$ be the moduli space of
$m$-truncated Barsotti--Tate groups in characteristic $p$ that
have codimension $c$ and dimension $d$. In other words, for each
$\FF_p$-scheme $S$, $\Bscr\Tscr_m(S)$ is the category of all $\BT_m$'s over
$S$ of codimension $c$ and dimension $d$, the morphisms in
$\Bscr\Tscr_m(S)$ being isomorphisms of $\BT_m$'s.

\noindent As explained in \cite{Wd},~Proposition~(1.8) and
Corollary~(3.3), it follows from results of Illusie and
Grothendieck (see~\cite{Ill}, Th\'eor\`eme~4.4) that $\Bscr\Tscr_m$ is  
a smooth
algebraic stack of finite type over $\FF_p$. More precisely,
$\Bscr\Tscr_m$ is an algebraic stack of the form
$[\pmb{GL}_{p^{m(c+d)}}\backslash Z_m]$, where $Z_m$ is a smooth
quasi-affine $\FF_p$-scheme on which $\pmb{GL}_{p^{m(c+d)}}$
acts. Moreover, the canonical morphism $P:Z_m \to \Bscr\Tscr_m$ is a
$\pmb{GL}_{p^{m(c+d)}}$-torsor for the Zariski topology. Thus $P(R):Z_m(R) \twoheadrightarrow
\Bscr\Tscr_m(R)$ is surjective for each commutative local $\FF_p$-algebra $R$.
\end{segment}

\begin{segment}{The level $m$ stratification}{LevelmStrata}
Let $k$ be a field of characteristic $p$ and let $D_m$ be a
$\BT_m$ over $k$. By the definition of the stack $\Bscr\Tscr_m$,
$D_m$ defines a
$1$-morphism over $k$
\[
\xi := \xi_{D_m}\colon \Spec k \to \Bscr\Tscr_m
\otimes_{\FF_p} k.
\]
The pair $(\xi,k)$ defines a point of $\Bscr\Tscr_m
\otimes_{\FF_p} k$ (in the sense of \cite{LMB}, Section (5.2))
which we also denote by $\xi$. As $\Bscr\Tscr_m$ is locally
noetherian, $\xi$ is algebraic by \cite{LMB}, Section (11.3)
and its residue field is $k$. Let
$\Gscr_{\xi}$ be the residue gerbe of the point $\xi$; it is an algebraic stack
which is an fppf gerbe over $\Spec k$.

\begin{proclamation}{Lemma}\label{LemResidue}
The canonical monomorphism $\Gscr_{\xi} \to \Bscr\Tscr_m\otimes_{\FF_p} k$ is  
representable by an
immersion of finite presentation. The algebraic stack $\Gscr_{\xi}$ is  
smooth over $\Spec k$.
\end{proclamation}

\begin{proof}
The morphism $\Gscr_{\xi} \to \Bscr\Tscr_m\otimes_{\FF_p} k$ is representable because  
$\xi$ is algebraic. For the remaining assertions we may assume that $k$ is algebraically closed. With the notations of Subsection \eqref{ModTruncBT}, from \cite{LMB}, Exemple (11.2.2) we get that the
fibre product of the diagram
\[\xymatrix{
& Z_m\otimes_{\dbF_p} k \ar[d]^{P_k} \\
\Gscr_{\xi} \ar[r] & \Bscr\Tscr_m\otimes_{\FF_p} k
}\]
is the $\pmb{GL}_{p^{m(c+d)}}$-orbit $O(x)$ in $Z_m\otimes_{\FF_p} k $, where $x \in Z_m(k)$ is a lift of $\xi$. As $O(x) \to Z_m\otimes_{\FF_p} k$ is a quasi-compact immersion of noetherian schemes, $\Gscr_{\xi} \to \Bscr\Tscr_m\otimes_{\FF_p} k$ is representable by an immersion of finite presentation.

The morphism $P_k$ is smooth and surjective and thus $O(x) \to
\Gscr_{\xi}$ is smooth and surjective. As $O(x)$ is smooth over
$k$, $\Gscr_{\xi}$ is smooth over $\Spec k$.
\end{proof}

\noindent Let $S$ be an arbitrary $k$-scheme. Let $X_m$ be a $\BT_m$ over
$S$ defining a \hbox{1-morphism} $\xi_{X_m}\colon S \to \Bscr\Tscr_m\otimes_{\FF_p} k$. Let
$S_{D_m}(X_m)$ be the fibre product of the diagram
\[\xymatrix{
& S \ar[d]^{\xi_{X_m}} \\
\Gscr_{\xi} \ar[r] & \Bscr\Tscr_m\otimes_{\FF_p} k. }\] The canonical morphism
$S_{D_m}(X_m) \to S$ is an immersion of finite presentation by
Lemma~\ref{LemResidue}. Thus we will view $S_{D_m}(X_m)$ as a
locally closed subscheme of $S$. As $S_{D_m}(X_m) \to S$ is quasi-compact, its
scheme-theoretic closure $\overline{S_{D_m}(X_m)}$ exists
by~\cite{EGAI}, Corollaire~(6.10.6).

By~\cite{LMB}, Section (11.1), a morphism $f\colon S' \to S$ of
$k$-schemes factors through $S_{D_m}(X_m)$ if and only if $f^*(X_m)$ and $D_m \times_{\Spec k} S'$ are locally for the fppf topology isomorphic as $\BT_m$'s over $S'$. We call $S_{D_m}(X_m)$ the \emph{level $m$ stratum of $(S,X_m)$ with respect
to $D_m$}.

\noindent We note that the level~$1$ strata are the Ekedahl--Oort
strata introduced in~\cite{Oo2} and that the level~$m$ strata were
studied first in~\cite{Wd} and in~\cite{Va2}.
\end{segment}

\begin{segment}{D\'evissage}{Dev}
We will show that it suffices to prove Theorem~\ref{MainThm} in
the Essential Situation~\ref{EssSit}. Generalizing
Definition~\ref{DefPure}, we say that a substack $\Tscr$ of an
algebraic stack $\Sscr$ is \emph{pure in $\Sscr$} if the immersion
$\Tscr \hookrightarrow \Sscr$ is affine. We recall the
following lemma (which follows from the fact that the affineness  
property for a morphism is local for the fpqc topology, see  
\cite{EGAIV}, Proposition
(2.7.1)).

\begin{proclamation}{Lemma}\label{PureProp}
Let $\Sscr$ be an algebraic stack and let $\Tscr \subset \Sscr$ be a  
substack. Let
$f\colon \Yscr \to \Sscr$ be a representable morphism of algebraic  
stacks. We have
the following two properties:
\medskip
\begin{assertionlist}
\item
If $\Tscr$ is pure in $\Sscr$, $f^{-1}(\Tscr)$ is pure in $\Yscr$.
\item
Conversely, assume that $f$ is quasi-compact and faithfully flat. If  
$f^{-1}(\Tscr)$
is pure in $\Yscr$, then $\Tscr$ is pure in $\Sscr$.
\end{assertionlist}
\end{proclamation}

We now refer to the general situation of
Subsection~\eqref{LevelmStrata}. It follows from the construction
of $S_{D_m}(X_m)$ that to prove that the immersion
$S_{D_m}(X_m)\hookrightarrow S$ is affine, it suffices to show that
the immersion $\Gscr_{\xi} \to \Bscr\Tscr_m$ is affine.

A smooth scheme $\Ascr$ over $k$ of dimension $cd$ with a
$p$-divisible group $\Escr$ over $\Ascr$ satisfies the properties
(a) and (b) of the Essential Situation \ref{EssSit} if and only if
the morphism $\xi_{\Escr[p^m]}\colon \Ascr \to \Bscr\Tscr_m
\otimes_{\FF_p} k$ defined by $\Escr[p^m]$ is smooth and contains
the image of $\xi_{D_m}$. Note that in this case, the level $m$
stratum $\grs_D(m)=\Ascr_D(\Escr[p^m])$ is smooth over $\Gscr_{\xi}$ and
hence by Lemma~\ref{LemResidue} over $k$. This was claimed in the
introduction.

To reduce to the Essential Situation~\ref{EssSit}, by  
Lemma~\ref{PureProp}, we can assume that $k$ is algebraically closed,  
and it suffices to prove the following proposition.

\begin{proclamation}{Proposition}\label{RedtoStand}
There exists a smooth $k$-scheme $\Ascr$ of finite type which is
equidimensional of dimension $cd$ and a $p$-divisible group
$\Escr$ of codimension $c$ and dimension $d$ over $\Ascr$ such
that for each $m\in\Bbb N^*$ the morphism $\xi_{\Escr[p^m]}\colon \Ascr \to \Bscr\Tscr_m\otimes_{\FF_p} k$ defined by $\Escr[p^m]$ is smooth and
surjective.
\end{proclamation}

\begin{proof}
For $\Ascr$, we will take the special fibre of a good integral
model of a Shimura variety $\text{Sh}(\Gscr,\Xscr)$ associated to
a certain PEL-datum as follows. If $c=d=1$, we can take
$\text{Sh}(\Gscr,\Xscr)$ to be the elliptic modular curve. Thus we
can assume that $r:=c+d\ge 3$; in this case the PEL-datum will be
unitary.

Let $\KK$ be a quadratic imaginary extension of $\QQ$ in which $p$ splits. Let $O_{\KK}$ be the ring of integers of
$\KK$. Let $\star$ be the nontrivial automorphism of $\KK$. Let
$\VV$ be a $\QQ$--vector space of dimension $2r$. We fix a
monomorphism $\Bbb K\hookrightarrow \End(\Bbb V)$ of $\Bbb Q$--algebras.
Via this monomorphism, we can view $\Bbb V$ naturally as a
$\KK$-vector space of dimension $r$ and we can view
$\text{Res}_{\KK/\Bbb Q} \Bbb G_{m,\Bbb K}$ as a torus of
$\pmb{GL}_{\Bbb V}$. Whenever we write $\pmb{SL}_{\Bbb V}$ or
$\pmb{GL}_{\Bbb V}$, we consider $\Bbb V$ only as a $\Bbb Q$--vector
space.

\noindent Let $\Gscr^{\text{der}}$ be the simply connected
semisimple group over $\Bbb Q$ whose $\Bbb Q$--valued points are those
$\Bbb K$-valued points of $\pmb{SL}_{\Bbb V}$ that leave invariant the
hermitian form
$-z_1z_1^{\star}-\cdots-z_cz_c^{\star}+z_{c+1}z_{c+1}^{\star}+\cdots
z_{r} z_{r}^{\star}$ on $\Bbb V$. The group
$\Gscr_{\Bbb R}^{\text{der}}$ is isomorphic to $\pmb{SU}(c,d)$.
Hilbert's Theorem 90 implies that there exists a unique (up to
non-zero scalar multiplication) symplectic form $\lrangle\colon
\VV \times \VV \to \QQ$ fixed by $\Gscr^{\text{der}}$. We have
$\langle bv,v'\rangle = \langle v,b\star v'\rangle$ for all $b \in
\Bbb K$ and $v,v' \in \VV$.

\noindent Let $\Gscr$ be the subgroup of $\pmb{GSp}(\Bbb V,<,>)$
generated by $\Gscr^{\text{der}}$ and by the torus
$\text{Res}_{\KK/\Bbb Q} \Bbb G_{m,\Bbb K}$. Our notations match i.e.,
$\Gscr^{\text{der}}$ is the derived group of $\Gscr$. It is easy
to see that there exists a $\Gscr(\RR)$-conjugacy class $\Xscr$ of
homomorphisms $h\colon \Res_{\CC/\RR}(\GG_{m,\CC}) \to
\Gscr_{\RR}$ such that every $h \in \Xscr$ defines a Hodge
$\Bbb Q$--structure on $\VV$ of type $\{(-1,0),(0,-1)\}$ (with the sign
convention of \cite{De}) and such that
\[
\VV_{\RR} \times \VV_{\RR} \to \RR, \qquad (v,v') \sends \langle  
v,h(\sqrt{-1})v'\rangle
\]
is symmetric and either positive or negative definite (for
instance, see \cite{De}, proof of Proposition 2.3.10 or \cite{Ko},
Lemma 4.3). Then $(\Gscr,\Xscr)$ is a Shimura pair given by the
PEL-datum $(\KK,\star,\VV,\lrangle,\Xscr)$. Its reflex field $E$
is either equal to $\QQ$ (if $c = d$) or isomorphic to $\KK$ (if
$c \ne d$). In both cases, for each prime $v$ of $E$ that divides
$p$ the completion $E_v$ of $E$ with respect to $v$ is $\QQ_p$.

\noindent As $p$ splits in $\Bbb K$, the reductive group
$\Gscr_{\Bbb Q_p}$ is split. This implies that there exists an
$O_{\KK}$-invariant $\ZZ_p$-lattice $\Gamma$ of $\VV \otimes_{\QQ}
\QQ_p$ such that the alternating form on $\Gamma$ induced by
$\lrangle$ is a perfect $\ZZ_p$-form. Denote by $\AA_f^{(p)}$ the
ring of finite ad\`eles of $\QQ$ with trivial $p$-th component and
fix an open compact subgroup $C^{(p)} \subset \Gscr(\AA_f^{(p)})$.
Let $\Mcal$ be the moduli space over $O_{E_v} = \ZZ_p$ of abelian
schemes associated to the data
$(\Bbb K,\star,\VV,\lrangle,O_{\Bbb K},\Gamma,C^{(p)})$, cf.
\cite{Ko}, Section 5.

\noindent We set $\Ascr := \Mcal \otimes_{\ZZ_p} k$. For $C^{(p)}$ small enough, $\Ascr$ is quasi-projective, smooth, and equidimensional of dimension $cd$ over $k$. As the
adjoint group of $\Gscr$ is simple, the Hasse principle holds for
$\Gscr$ (cf. \cite{Ko}, top of p. 394) and this implies that
$\Mcal$ is an integral model over $\Bbb Z_p$ of
$\text{Sh}(\Gscr,\Xscr)$ alone (cf. \cite{Ko}, Section 7).

\noindent Let $\Zscr$ be the $p$-divisible group of the universal
abelian scheme over $\Ascr$. The height of $\Zscr$ is $2r$ and its
dimension is $r$. It is endowed with an action of $O_{\KK}
\otimes_{\ZZ} \ZZ_p = \ZZ_p \times \ZZ_p$ and with a principal
quasi-polarization $\Lambda_{\Zscr}$. The action of $\ZZ_p \times
\ZZ_p$ defines a decomposition $\Zscr = \Escr \times \Escr'$,
where $\Escr$ is a $p$-divisible group over $\Ascr$ of codimension
$c$ and dimension $d$ and where $\Escr'$ is via
$\Lambda_{\Zscr}$ isomorphic to the Cartier dual $\Escr\vdual$ of $\Escr$. In
particular $\Zscr \cong \Escr \times \Escr\vdual$, endowed with
its natural $\ZZ_p \times \ZZ_p$-action and with its natural
principal quasi-polarization.

\noindent Due to the moduli interpretation of $\Mscr$, we easily
get that:

\medskip
{\bf (*)} if $y\colon \Spec k\to \Ascr$ and if $\tilde E_y$ is a  
$p$-divisible group isogenous to \hbox{$E_y:=y^*(\Escr)$}, then there  
exists a point $\tilde y:\Spec k\to \Ascr$ such that $\tilde E_y$ is  
isomorphic to $E_{\tilde y}:=\tilde y^*(\Escr)$.

\medskip
\noindent We claim that the morphism
\hbox{$\xi_{\Escr[p^m]}\colon \Ascr \to \Bscr\Tscr_m \otimes_{\FF_p} k$}
is smooth and surjective. To check this, let $R$ be a
local, artinian $k$-algebra and let $I \subset R$ be an ideal with
$I^2 = 0$. We define $R_0 := R/I$. Assume that we are given a
commutative diagram
\[\xymatrix{
\Spec R_0 \ar[r]^{\chi_0} \ar[d] & \Ascr \ar[d]^{\xi_{\Escr[p^m]}} \\
\Spec R \ar[r]^-{\upsilon} & \Bscr\Tscr_m \otimes_{\FF_p} k. }\]
By a theorem of Grothendieck (see \cite{Ill}, Th\'eor\`eme 4.4)
there exists a\break \hbox{$p$-divisible} group $E$ over $R$ which
lifts $\chi_0^*(\Escr)$ and such that $E[p^m]$ is the $\BT_m$
corresponding to $\upsilon$. We endow $E \times E\vdual$ with the
natural $\ZZ_p \times \ZZ_p$-action and with its natural principal
quasi-polarization. As explained above, we have $(E \times
E\vdual)\times_R R_0 \cong \chi_0^*(\Zscr)$. From this and the
Serre--Tate deformation theory, we get that there exists an
abelian scheme over $R$ whose $p$-divisible group is $E \times
E\vdual$ and such that its reduction modulo $I$ is given by
$\chi_0$, and this, due to the moduli interpretation of $\Mscr$, defines a morphism $\chi\colon \Spec R \to\Ascr$ that lifts $\chi_0$. Thus $\xi_{\Escr[p^m]}$ is
smooth.

By the generalization of the integral Manin principle to
certain Shimura varieties of PEL-type proved in \cite{Va4},
Subsection~5.4 (in particular Example~5.4.3~(b)),
$\xi_{\Escr[p^m]}$ is also surjective. For the reader's
convenience and as our present context is much simpler than the
general situation considered in \cite{Va4}, we give a direct proof
of the surjectivity of $\xi_{\Escr[p^m]}$. 

Let $K$ be an
algebraically closed extension of $k$ and let $\tilde\xi\colon \Spec K
\to \Bscr\Tscr_m\otimes_{\Bbb F_p} k$ be a $K$-valued point corresponding to a $\BT_m$ $\tilde D_m$ over $K$. Let $\tilde D$ be a
$p$-divisible group over $K$ such that $\tilde D[p^m] = \tilde D_m$.  
To show that $\xi_{\Escr[p^m]}$ is surjective it suffices to prove  
that there exists a point $y_3\colon \Spec K \to \Ascr$ such that the  
$p$-divisible groups $y_3^*(\Escr)$ and $\tilde D$ are isomorphic (thus the fibre product of $\tilde\xi$ and $\xi_{\Escr[p^m]}$ is non-empty). We check this in four steps as follows. 

{\bf (i)} Let $\Tscr_0$ be a maximal torus of $\Gscr^{\text{der}}$  
such that $\Tscr_{0,\Bbb R}$ is compact and $\Tscr_{0,\Bbb Q_p}$ has  
$\Bbb Q_p$-rank $0$, cf. \cite{Ha}, Lemma 5.5.3. Let $\Tscr$ be the  
unique maximal torus of $\Gscr$ that contains $\Tscr_0$ (it is  
generated by $\Tscr_0$ and by the center of $\Gscr$). Let  
$h_0\in\Xscr$ be such that it factors through $\Tscr_0$ (it exists as  
all maximal compact tori of $\Gscr^{\text{der}}_{\Bbb R}$ are  
$\Gscr^{\text{der}}(\Bbb R)$-conjugate and as the centralizer of each  
$h\in\Xscr$ in $\Gscr_{\Bbb R}^{\text{der}}$ is a maximal compact, connected subgroup of $\Gscr_{\Bbb R}^{\text{der}}$). We have an injective  
map $(\Tscr,\{h_0\})\hookrightarrow (\Gscr,\Xscr)$ of Shimura pairs.  
Each point of $\Mscr_{E_v}$ which is in the image of the natural  
functorial morphism $\text{Sh}(\Tscr,\{h_0\})_{E_v}\to\Mscr_{E_v}$  
specializes to a point $y_0:\Spec k\to\Ascr$ such that $y_0^*(\Escr)$  
is isoclinic (due to the fact that $\Tscr_{0,\Bbb Q_p}$ has  
$\Bbb Q_p$-rank $0$). Such a specialization makes sense as abelian  
varieties with complex multiplication over number fields have  
potentially good reduction everywhere.

{\bf (ii)} Based on (*) (applied over $K$) and (i) we get that there  
exists a point $y_1\colon \Spec K \to \Ascr$ such that $y_1^*(\Escr)$  
is isoclinic and its $a$-number is $1$.

{\bf (iii)} Based on (ii) and Grothendieck's specialization conjecture  
for $p$-divisible groups over $K$ of $a$-number $1$ (proved in  
\cite{Tr2}, Sections~6,~7, and~24; see also \cite{Oo1}, Theorem 6.2),  
there exists a point $y_2\colon \Spec K\to\Ascr$ such that the Newton  
polygons of $y_2^*(\Escr)$ and $\tilde D$ coincide. Here we are using the fact that $\xi_{\Escr[p^m]}$ is smooth.

{\bf (iv)} Based on (iii) and (*) there exists a point $y_3\colon  
\Spec K \to \Ascr$ such that the $p$-divisible groups $y_3^*(\Escr)$  
and $\tilde D$ are isomorphic.
\end{proof}
\end{segment}


\section{Group actions}\label{GroupActions}

From now on, we will be in the Essential Situation~\ref{EssSit}.
Thus $k$ is algebraically closed and $D$ is a $p$-divisible over
$k$ of codimension $c$ and dimension $d$. The height of $D$ is
$r:=c+d$. In this Section, we recall from \cite{Va2} the
definition of an action $\Bbb T_m\colon \Hscr_m \times \Dscr_m \to \Dscr_m$ of a linear algebraic group $\Hscr_m$ over $k$ on a $k$-scheme $\Dscr_m$ whose orbits parameterize isomorphism classes of $\BT_m$'s over $k$.

For a commutative $\Bbb F_p$-algebra $R$, let $W_m(R)$ be the ring
of Witt vectors of length $m$ with coefficients in $R$, let $W(R)$
be the ring of Witt vectors with coefficients in $R$, and let
$\sigma_R$ be the Frobenius endomorphism of either $W_m(R)$ or
$W(R)$ induced by the Frobenius endomorphism $r \mapsto r^p$ of
$R$. We set $\sigma:=\sigma_k$. Let $B(k)$ be the field of fractions  
of $W(k)$.

Let $(M,\phi)$ be the contravariant Dieudonn\'e module of $D$.
Thus $M$ is a free $W(k)$-module of rank $r$ and $\phi:M\to M$ is
a $\sigma$-linear endomorphism such that we have $pM\subseteq
\phi(M)$. Let $\vartheta:=p\phi^{-1}:M\to M$ be the Verschiebung
map of $(M,\phi)$. Let $M=F^1\oplus F^0$ be a direct sum
decomposition such that $\bar F^1:=F^1/pF^1$ is the kernel of the
reduction modulo $p$ of $\phi$. Let $\bar F^0:=F^0/pF^0$. The
ranks of $F^1$ and $F^0$ are $d$ and $c$ (respectively). We have
$\phi({1\over p}F^1\oplus F^0)=M$. The decomposition $M=F^1\oplus
F^0$ gives birth naturally to a direct sum decomposition of
$W(k)$-modules
$$\End(M)=\Hom(F^0,F^1)\oplus\End(F^1)\oplus\End(F^0)\oplus\Hom(F^1,F^0).$$
The association $e\to \phi(e):=\phi\circ e\circ \phi^{-1}$ defines a  
$\sigma$-linear
automorphism $\phi:\End(M)[{1\over p}]\arrowsim \End(M)[{1\over p}]$  
of $B(k)$-algebras.

\noindent Let $\Wscr_+$ be the maximal subgroup scheme of
$\pmb{GL}_{M}$ that fixes both $F^1$ and $M/F^1$; it is a closed
subgroup scheme of $\pmb{GL}_M$ whose Lie algebra is the direct
summand $\Hom(F^0,F^1)$ of $\End(M)$ and whose relative dimension
is $cd$. Let $\Wscr_0:=\pmb{GL}_{F^1}\times_{W(k)}
\pmb{GL}_{F^0}$; it is a closed subgroup scheme of $\pmb{GL}_M$
whose Lie algebra is the direct summand $\End(F^1)\oplus\End(F^0)$
of $\End(M)$ and whose relative dimension is $d^2+c^2$. The
maximal parabolic subgroup scheme $\Wscr_{+0}$ of $\pmb{GL}_{M}$
that normalizes $F^1$ is the semidirect product of $\Wscr_+$ and
$\Wscr_0$. Let $\Wscr_-$ be the maximal subgroup scheme of
$\pmb{GL}_{M}$ that fixes $F^0$ and $M/F^0$; it is a closed
subgroup scheme of $\pmb{GL}_M$ whose Lie algebra is the direct
summand $\Hom(F^1,F^0)$ of $\End(M)$ and whose relative dimension
is $cd$. The maximal parabolic subgroup scheme $\Wscr_{0-}$ of
$\pmb{GL}_{M}$ that normalizes $F^0$ is the semidirect product of
$\Wscr_-$ and $\Wscr_0$. If $R$ is a commutative $W(k)$-algebra, then we have:
$$\Wscr_+(R)=1_{M\otimes_{W(k)} R}+\Hom(F^0,F^1)\otimes_{W(k)} R$$
$$\Wscr_-(R)=1_{M\otimes_{W(k)} R}+\Hom(F^1,F^0)\otimes_{W(k)} R.$$
These identities imply that the group schemes $\Wscr_+$ and $\Wscr_-$  
are isomorphic
to $\Bbb G_a^{cd}$ over $\Spec W(k)$; in particular, they are smooth and  
commutative.
Let
$$\Hscr:=\Wscr_+\times_{W(k)}\Wscr_0\times_{W(k)} \Wscr_-;$$
it is a smooth, affine scheme  over $\Spec W(k)$ of relative dimension
$cd+d^2+c^2+cd=r^2$. We consider the natural product morphism
$\Pscr_0:\Hscr\to\pmb{GL}_{M}$ and the following morphism
$\Pscr_-:=1_{\Wscr_+}\times 1_{\Wscr_0}\times p1_{\Wscr_-}:\Hscr\to\Hscr$. Let
$$\Pscr_{0-}:=\Pscr_0\circ\Pscr_{-}:\Hscr\to\pmb{GL}_M;$$
it is a morphism of $\Spec W(k)$-schemes whose generic fibre is an  
open embedding of
$\Spec B(k)$-schemes.

\noindent Let $\tilde\Hscr$ be the dilatation of
$\pmb{\text{GL}}_M$ centered on $\Wscr_{+0,k}$ (see \cite{BLR},
Chapter~3, Section 3.2 for dilatations). We recall that if
$\pmb{GL}_M=\Spec R_M$ and if $I_{+0k}$ is the ideal of $R_M$ that
defines $\Wscr_{+0,k}$, then as a scheme $\tilde\Hscr$ is the
spectrum of the $R_M$-subalgebra $R_{\tilde\Hscr}$ of $R_M[{1\over
p}]$ generated by all elements ${*\over p}$ with $*\in I_{+0k}$.
It is well known that $\tilde\Hscr$ is a smooth,
affine group scheme over $\Spec W(k)$ which is uniquely determined
by the following two additional properties (they follow directly
from the definition of $R_{\tilde\Hscr}$; see Propositions~1, 2,
and 3 of loc.~cit.):

\medskip\begin{romanlist}
\item
There exists a homomorphism $\tilde\Pscr_{0-}:\tilde\Hscr\to \pmb{GL}_M$ whose
generic fibre is an isomorphism of $\Spec B(k)$-schemes.
\item
A morphism $f:Y\to \pmb{GL}_M$ of flat $\Spec W(k)$-schemes factors (uniquely)
through $\tilde\Pscr_{0-}$ if and only if the morphism $f_k:Y_k\to  
\pmb{GL}_{M/pM}$
factors through $\Wscr_{+0,k}$.
\end{romanlist}

\medskip

\noindent The group $\tilde\Hscr(W(k))$ is the parahoric subgroup
of $\pmb{\text{GL}}_M(W(k))$ that normalizes the sublattice $F_1
\oplus pF_0=F^1+pM$ of $M$.

\interbreak

\noindent The morphism $\Pscr_{0-}$ factors naturally as a
morphism $\Pscr:\Hscr\to\tilde\Hscr$ (cf. (ii)) whose $p$-adic completion is
an isomorphism (cf. \cite{Va2}, Subsubsection 2.1.1). Therefore we
have a natural identification
$\Hscr_{W_m(k)}=\tilde\Hscr_{W_m(k)}$ that provides
$\Hscr_{W_m(k)}$ with a group scheme structure over $\Spec W_m(k)$ which
does not depend on the decomposition
$\Hscr=\Wscr_+\times_{W(k)}\Wscr_0\times_{W(k)} \Wscr_-$ produced
by the choice of the direct sum decomposition $M=F^1\oplus F^0$.

\noindent For $g\in\pmb{GL}_M(W(k))$ and
$h=(h_1,h_2,h_3)\in\Hscr(W(k))$, let $g[m]\in\pmb{GL}_M(W_m(k))$
and $h[m]=(h_1[m],h_2[m],h_3[m])\in\Hscr(W_m(k))$ be the
reductions  modulo $p^m$ of $g$ and $h$ (respectively). Thus
$1_{M/p^mM}=1_M[m]$. Let $\phi_m,\vartheta_m:M/p^mM\to M/p^mM$ be
the reductions  modulo $p^m$ of $\phi,\vartheta:M\to M$.

\noindent Let $\sigma_{\phi}:M\arrowsim M$ be the $\sigma$-linear
automorphism which takes $x\in F^1$ to ${1\over p}\phi(x)$ and
takes $x\in F^0$ to $\phi(x)$. Let $\sigma_{\phi}$ act on the sets
underlying the groups $\pmb{GL}_M(W(k))$ and $\pmb{GL}_M(W_m(k))$
in the natural way: if $g\in\pmb{GL}_M(W(k))$, then
$\sigma_{\phi}(g):=\sigma_{\phi} g\sigma_{\phi}^{-1}$ and
$\sigma_{\phi}(g[m]):=(\sigma_{\phi} g\sigma_{\phi}^{-1})[m]$. For
$g\in\Wscr_+(W(k))$ (resp. $g\in\Wscr_0(W(k))$ or
$g\in\Wscr_-(W(k))$) we have $\phi(g)=\sigma_{\phi}(g^p)$ (resp.
we have $\phi(g)=\sigma_{\phi}(g)$ or
$\phi(g^p)=\sigma_{\phi}(g)$).

Let $\text{Aff}_k$ be the category of affine schemes over $k$. Let
$\text{Set}$ and $\text{Group}$ be the categories of abstract sets
and groups (respectively). Let $\vartriangle$ be a smooth, affine (resp. a
smooth, affine group) scheme of finite over $\Spec W(k)$. Let
$\Bbb W_m(\vartriangle):\text{Aff}_k\to\text{Set}$ (resp.
$\Bbb W_m(\vartriangle):\text{Aff}_k\to\text{Group}$) be the contravariant
functor that associates to an affine $k$-scheme $\Spec R$ the set
(resp. the group) $\vartriangle(W_m(R))$. This functor is representable by
an affine, smooth (resp. affine, smooth group) scheme over $k$ of
finite type to be denoted also by $\Bbb W_m(\vartriangle)$ (see \cite{Va2},
Subsection 2.1.4 for these facts due to Greenberg).

\noindent
Let $\Hscr_m:=\Bbb W_m(\Hscr)$ and $\Dscr_m:=\Bbb W_m(\pmb{GL}_{M})$.
As $\Pscr_{W_m(k)}:\Hscr_{W_m(k)}\to\tilde\Hscr_{W_m(k)}$
is an isomorphism of $\Spec W_m(k)$-schemes, we will identify
naturally
$$\Hscr(W(k))=\tilde\Hscr(W(k))\;\;\text{and}\;\;\Hscr_m=\Bbb W_m(\Hscr)=\Bbb W_m(\tilde\Hscr).$$
Thus in what follows we will view $\Hscr(W(k))$ as a subgroup of  
$\pmb{GL}_M(W(k))$
and $\Hscr_m$ as a connected, smooth, affine group over $k$ of  
dimension $mr^2$ (cf.
\cite{Va2}, Subsection 2.1.4 applied to $\tilde\Hscr$). Similarly, we  
will view
$\Dscr_m$ as a connected, smooth, affine variety over $k$ of dimension  
$mr^2$. Let
\[
\Bbb T_m\colon \Hscr_m\times_k \Dscr_m\to \Dscr_m
\]
be the action defined on $k$-valued points as follows. If  
$h=(h_1,h_2,h_3)\in\Hscr(W(k))$ and
$g\in \pmb{GL}_{M}(W(k))$, then the product of $h[m]=(h_1[m],h_2[m],h_3[m])\in
\Hscr_m(k)=\Hscr(W_m(k))$ and $g[m]\in\Dscr_m(k)=\pmb{GL}_M(W_m(k))$  
is the element
$$\Bbb T_m(h[m],g[m]):=(h_1h_2h_3^pg\phi(h_1h_2h_3^p)^{-1})[m]$$
$$=(h_1h_2h_3^pg\phi(h_3^p)^{-1}\phi(h_2)^{-1}\phi(h_1)^{-1})[m]=(h_1h_2h_3^pg\sigma_{\phi}(h_3)^{-1}\sigma_{\phi}(h_2)^{-1}\sigma_{\phi}(h_1^p)^{-1})[m]$$
$$=h_1[m]h_2[m]h_3[m]^pg[m]\sigma_{\phi}(h_3[m])^{-1}\sigma_{\phi}(h_2[m])^{-1}\sigma_{\phi}(h_1[m]^p)^{-1}\in\Dscr_m(k).$$
The formula
$\Bbb T_m(h[m],g[m])=(h_1h_2h_3^pg\phi(h_1h_2h_3^p)^{-1})[m]$ shows
that the action $\Bbb T_m$ is intrinsically associated to $D$ i.e.,
it does not depend on the choice of the direct sum decomposition
$M=F^1\oplus F^0$.

Let $\Oscr_m$ be the orbit of $1_M[m]\in\Dscr_m(k)$ under the
action $\Bbb T_m$. Let $\bar\Oscr_m$ be the scheme-theoretic closure of
$\Oscr_m$ in $\Dscr_m$; it is an affine, integral scheme over $k$.
The orbit $\Oscr_m$ is a connected, smooth, open subscheme of
$\bar\Oscr_m$ and thus it is also a quasi-affine scheme over
$k$. Let $\Sscr_m$ be the subgroup scheme of $\Hscr_m$ which is
the stabilizer of $1_M[m]$ under the action $\Bbb T_m$. Let
$\Cscr_m$ be the reduced group of $\Sscr_m$. Let $\Cscr_m^0$ be
the identity component of $\Cscr_m$. The connected, smooth group
$\Cscr_m^0$ over $k$ is {\it unipotent} i.e., it has no torus of positive dimension (see \cite{Va2},
Theorem 2.4 (a); see also Subsection~\eqref{CompC10}).

\section{Combinatorics of the action $\Bbb T_1$}\label{Combinatorics}

In this Section, we present basic combinatorial properties of the action
$\Bbb T_1$ which will be used in Section \ref{MainProof} to show that if the condition $\Cond$
holds for $D[p]$, then the orbit $\Oscr_m$ of $\Bbb T_m$ is affine. For the
remaining part of this Section,we let $m=1$ and we use the notations of Section \ref{GroupActions}.

\begin{segment}{Nilpotent subalgebras of $\End(M/pM)$}{Nilp}
In \cite{Kr} (see also \cite{Oo2}, Subsection (2.3) and Lemma
(2.4)) it is shown that there exists a $k$-basis $\{\bar
e_1,\ldots,\bar e_r\}$ for $M/pM$ and a permutation $\pi$ of the
set $J:=\{1,\ldots,r\}$ such that for each $i\in J$, the
following two properties hold:

\medskip\begin{romanlist}
\item
$\phi_1(\bar e_i)=0$ if $i>c$, and $\phi_1(\bar e_i)=\bar e_{\pi(i)}$ if $i\le c$;

\item
$\vartheta_1(\bar e_{\pi(i)})=0$ if $i\le c$, and $\vartheta_1(\bar e_{\pi(i)})=\bar e_i$ if $i>c$.

\end{romanlist}

\medskip
\noindent The permutation $\pi$ is not uniquely
determined by the isomorphism class of $D[p]$. For instance, we
can always replace $\pi$ by $\pi_0\pi_1\pi\pi_0^{-1}$, where
$\pi_0$ is an arbitrary permutation of the set $J$ that leaves
invariant the subset $\{1,\ldots,c\}$ and where $\pi_1$ is an
arbitrary permutation of the set $J$ that leaves invariant all the
subsets $\{\pi^i(1),\ldots,\pi^i(c)\}$ with
$i\in\Bbb N$. More precisely, it is known
(e.g.,~\cite{Mo},~\cite{Va3}, or~\cite{MW}) that there exists a
canonical bijection between isomorphism classes of $\BT_1$'s over $k$  
and the quotient set $S_r/(S_c\times S_d)$.

\noindent Let $\{e_1,\ldots,e_r\}$ be a $W(k)$-basis for $M$ that
lifts the $k$-basis $\{\bar e_1,\ldots,\bar e_r\}$ for $M/pM$ and
such that $F^1=\oplus_{i=c+1}^r W(k)e_i$. Let $\{e_{i,j}|i,j\in
J\}$ be the $W(k)$-basis for $\End(M)$ such that for each $l\in J$ we
have $e_{i,j}(e_l)=\delta_{j,l}e_i$. Let $\{\bar e_{i,j}|i,j\in
J\}$ be the the reduction
modulo $p$ of $\{e_{i,j}|i,j\in J\}$; it is a $k$-basis for $\End(M/pM)$. Let $\sigma_{\pi}:M\arrowsim
M$ be the $\sigma$-linear automorphism that maps $e_i$ to
$e_{\pi(i)}$ for all $i\in J$. Let
$g_{\pi}:=\sigma_{\pi}\sigma_{\phi}^{-1}\in\pmb{GL}_M(W(k))$. Due
to the properties (i) and (ii), the reduction modulo $p$ of
$(M,g_{\pi}\phi,\vartheta g_{\pi}^{-1})$ coincides with
$(M/pM,\phi_1,\vartheta_1)$. Based on this, we can assume that
$g_{\pi}[1]=1_M[1]$; thus $\sigma_{\phi}$ and $\sigma_{\pi}$ are
congruent modulo $p$. As the action $\Bbb T_1$ is intrinsically
associated to $D$ (i.e., it does not depend on the choice of the
direct sum decomposition $M=F^1\oplus F^0$), to study the group
$\Cscr_1^0$ we can assume $F^0=\oplus_{i=1}^c W(k)e_i$. Let
$$\Jscr_+:=\{(i,j)\in J^2|j\le c<i\},\;\; $$
$$\Jscr_0:=\{(i,j)\in J^2|\;\text{either}\;i,j>c\;\text{or}\;i,j\le
c\},\;\;\;\text{and}\;\;\;\Jscr_-:=\{(i,j)\in J^2|i\le c<j\}.$$
\noindent The three sets $\{\bar e_{i,j}|(i,j)\in\Jscr_+\}$,
$\{\bar e_{i,j}|(i,j)\in\Jscr_0\}$, and $\{\bar
e_{i,j}|(i,j)\in\Jscr_-\}$ are $k$-bases for
$\Lie(\Wscr_{+,k})=\Hom(\bar F^0,\bar F^1)$,
$\Lie(\Wscr_{0,k})=\End(\bar F^1)\oplus\End(\bar F^0)$, and
$\Lie(\Wscr_{-,k})=\Hom(\bar F^1,\bar F^0)$ (respectively).

\noindent By the {\it $\pi$-order} of the pair $(i,j)\in\Jscr_-$,
we mean the smallest positive integer $\nu(i,j)$ such that we have
$$(\pi^{\nu(i,j)}(i),\pi^{\nu(i,j)}(j))\in\Jscr_+\cup\Jscr_-.$$
We define the following five sets:
$$\Jscr_{-,1}:=\{(i,j)\in\Jscr_-|(\pi^{\nu(i,j)}(i),\pi^{\nu(i,j)}(j))\in\Jscr_+\}\;\;\;\text{and}\;\;\;\Jscr_{-,2}:=\Jscr_-\setminus\Jscr_{-,1}$$
$$\Jscr_{+,1}:=\{(\pi^{\nu(i,j)}(i),\pi^{\nu(i,j)}(j))|(i,j)\in\Jscr_{-,1}\}\;\;\;\text{and}\;\;\;\Jscr_{+,2}:=\Jscr_+\setminus\Jscr_{+,1}$$
$$\Jscr_{0,0}:=\{(\pi^s(i),\pi^s(j)|(i,j)\in\Jscr_{-,1}\;\;\text{and}\;\;s\in\{1,
\dots, \nu(i,j)-1\}\}.$$
We remark that the set $\Jscr_+\cup\Jscr_{0,0}\cup\Jscr_-$ contains no  
pair of the
form $(i,i)$. The number of elements of the set $\Jscr_{0,0}$ is
$$|\Jscr_{0,0}|:=\sum_{(i,j)\in\Jscr_{-,1}} \Big( \nu(i,j)-1 \Big).$$

\noindent For $(i,j)\in\Jscr_{-,1}$ and
$s\in\{0,1,\ldots,\nu(i,j)\}$ we define the {\it
$\pi$-level} of $(\pi^s(i),\pi^s(j))$ to be the number
$\eta(\pi^s(i),\pi^s(j)):=s$ and we define the $\pi$-order of
$(\pi^s(i),\pi^s(j))$ to be the number
$\nu(\pi^s(i),\pi^s(j)):=\nu(i,j)-s$. Thus the $\pi$-order
$\nu(i,j)$ and the $\pi$-level $\eta(i,j)$ are well defined for
all pairs $(i,j)\in\Jscr_{+,1}\cup\Jscr_{0,0}\cup\Jscr_{-,1}$.

\begin{proclamation}{Lemma}\label{NilLem}
Let $\nfr_{\bigstar}$ be the $k$-vector space generated by those $\bar  
e_{(i,j)}$'s
with $(i,j)\in \Jscr_{\bigstar}$ (thus  
$\nfr_{+,1}:=\oplus_{(i,j)\in\Jscr_{+,1}}
k\bar e_{i,j}$, $\nfr_{0,0}:=\oplus_{(i,j)\in\Jscr_{0,0}}
k\bar e_{i,j}$, etc.). Then the following five properties hold:
\begin{definitionlist}
\item If the pairs $(i,j)$ and $(j,l)$ belong to $\Jscr_{0,0}$,
then we have $(i,l)\in\Jscr_{0,0}$. Similarly, if one of the pairs
$(i,j)$ and $(j,l)$ belongs to $\Jscr_{0,0}$ and the other one
belongs to $\Jscr_{+,1}$ (resp. to $\Jscr_{-,1}$), then  we have
$(i,l)\in\Jscr_{+,1}$ (resp. we have $(i,l)\in\Jscr_{-,1}$). \item
The $k$-vector space $\nfr_{0,0}$ is a nilpotent subalgebra of
$\End(M/pM)$. More precisely, we have
$\nfr_{0,0}^{\text{max}(c,d)}=0$. \item The $k$-vector spaces
$\nfr_+$, $\nfr_{+,1}$, $\nfr_-$, and $\nfr_{-,1}$ are both left
and right $\nfr_{0,0}$-modules. Moreover we have
$\nfr_{0,0}^d\nfr_+=\nfr_+\nfr_{0,0}^c=\nfr_{0,0}^c\nfr_-=\nfr_-\nfr_{0,0}^d=0$.
\item The $k$-vector spaces $\nfr_+\oplus\nfr_{0,0}$,
$\nfr_{+,1}\oplus\nfr_{0,0}$, $\nfr_{0,0}\oplus\nfr_-$, and
$\nfr_{0,0}\oplus\nfr_{-,1}$ are nilpotent subalgebras of
$\End(M/pM)$.
\end{definitionlist}
\end{proclamation}

\begin{proof}
We prove only the first part of (a) as the second part of (a)
is proved similarly. Let
$s:=\min\{\eta(i,j),\eta(j,l)\}\in\Bbb N^*$ and let
$t:=\min\{\nu(i,j),\nu(j,l)\}\in\Bbb N^*$. From the very definition of
$s$, we get that one of the two pairs $(\pi^{-s}(i),\pi^{-s}(j))$ and
$(\pi^{-s}(j),\pi^{-s}(l))$ belongs to $\Jscr_-$ while the other
one belongs to $\Jscr_0$; thus $(\pi^{-s}(i),\pi^{-s}(l))\in\Jscr_-$. It is easy to see that for all
$u\in\{1,\ldots,s+t-1\}$ we
have $(\pi^{-s+u}(i),\pi^{-s+u}(l))\in\Jscr_0$. From the very  
definition of $t$, we get
that one of the two pairs $(\pi^{t}(i),\pi^{t}(j))$ and
$(\pi^{t}(j),\pi^{t}(l))$ belongs to $\Jscr_+$ while the other one
belongs to $\Jscr_0$; thus $(\pi^{t}(i),\pi^{t}(l))\in\Jscr_+$. From the last three sentences
we get that we have $(\pi^{-s}(i),\pi^{-s}(l))\in\Jscr_{-,1}$ and
$\nu(\pi^{-s}(i),\pi^{-s}(l))=s+t$. As $1\le s<s+t$, we have
$(i,l)\in\Jscr_{0,0}$ (i.e., the first part of (a) holds) as well
as $\eta(i,l)=s$ and $\nu(i,l)=t$.

We check (b). For $\nu\in\Bbb N$ let $\nfr_{0,0}^\nu$ be the  
$k$-span of those $\bar e_{i,j}$'s for which $(i,j) \in \Jscr_{0,0}$  
and $\nu(i,j) = \nu$.
As $\nfr_+\nfr_+=0$, we have $\nfr_{0,0}^\nu \nfr_{0,0}^\nu = 0$. From the end of the  
previous paragraph we get that for $\nu_1 <
\nu_2$, we have: $\nfr_{0,0}^{\nu_1} \nfr_{0,0}^{\nu_2} +
\nfr_{0,0}^{\nu_2} \nfr_{0,0}^{\nu_1} \subseteq\nfr_{0,0}^{\nu_1}$.  
These imply that $\nfr_{0,0}$ is a nilpotent subalgebra
of $\End(M/pM)$. Thus $\nfr_{0,0}$ is the Lie algebra
of the smooth, connected subgroup of $\Wscr_{0,k}=\pmb{GL}_{\bar
F^1}\times_k \pmb{GL}_{\bar F^0}$ whose group of valued points in an  
arbitrary commutative $k$-algebra $R$ is $1_{M\otimes_{W(k)}  
R}+\nfr_{0,0}\otimes_k R$. Therefore $\nfr_{0,0}$ is
$\Wscr_0(k)$-conjugate to a Lie subalgebra of $\End(\bar
F^1)\oplus\End(\bar F^0)$ formed by upper triangular nilpotent
matrices. This implies that $\nfr_{0,0}^{\text{max}(c,d)}=0$. Thus
(b) holds.

It is obvious that $\nfr_+$ and $\nfr_-$ are left and right  
$\nfr_{0,0}$-modules. As in the previous paragraph, using the second  
part of (a) we argue that
$\nfr_{+,1}$ and $\nfr_{-,1}$ are left and right $\nfr_{0,0}$-modules.  
The second
part of (c) follows from relations of the form $\nfr_{0,0}^d\nfr_+\subseteq
[\nfr_{0,0}\cap \End(\bar F^1)]^d\Hom(\bar F^0,\bar F^1)=0$. Thus (c) holds.

Part (d) follows from (b) and (c) and the fact that $\nfr_+^2=\nfr_-^2=0$.
\end{proof}

\end{segment}

\begin{segment}{Condition \Cond}{COND}
We define four subsets of $\cup_{s\ge 2} J^s$ as follows:
$$\Gamma:=\{(i_1,\ldots,i_s)|s\ge  
2,(i_{\ell},i_{\ell+1})\in\Jscr_{0,0} \  \forall \
\ell\in\{1,\ldots,s-2\}, (i_{s-1},i_s)\in\Jscr_{+,2}\}$$
$$\Delta:=\{(i_1,\ldots,i_s)|s\ge
3,(i_1,\ldots,i_{s-1})\in\Gamma,(i_{s-1},i_s)\in\Jscr_{0,0}\}$$
$$\Gamma_1:=\{(i_1,\ldots,i_s)\in\Gamma|(i_1,i_s)\in\Jscr_{+,1}\;\text{and}\;(i_2,i_s)\notin\Jscr_{+,1}\}$$
$$\Delta_1:=\{(i_1,\ldots,i_s)\in\Delta|(i_1,i_s)\in\Jscr_{+,1}\;\text{and}\;(i_2,i_s)\notin\Jscr_{+,1}\}.$$
If $(i_1,\ldots,i_s)\in\Gamma_1$, then $s\ge 3$. For each element $\gamma=(i_1,\ldots,i_s)\in\Gamma\cup\Delta$ and for every $t\in \Bbb N^*$, let $n_t(\gamma)$ be the number of
elements $\ell\in\{1,\ldots, s-1\}$ such that we have
$(i_{\ell},i_{\ell+1})\in\Jscr_{0,0}$ and
$\nu(i_{\ell},i_{\ell+1})=t$. We define a number:
$$\kappa(\gamma):=\sum_{t=1}^{\infty} n_t(\gamma)p^{-t}\in  
(0,\infty)\cap\Bbb Q.$$

\begin{proclamation}{Definition}\label{Defkappa}
Let $\kappa(\pi):=\max\{\kappa(\gamma)|\gamma\in
\Gamma_1\cup\Delta_1\}$ with the convention that $\kappa(\pi)=0$ if $\Gamma_1\cup\Delta_1$ is the empty set. Let $\kappa(D[p])$ be the smallest value
$\kappa(\pi)$, where $\pi$ runs through all permutations of the
set $\{1,\ldots,r\}$ for which there exists a $k$-basis $\{\bar
e_1,\ldots,\bar e_r\}$ for $M/pM$ as in the beginning of
Subsection \eqref{Nilp}. We say the condition $\Cond$ holds for
$D$ (or for any truncation of $D$) if either $\kappa(D[p])<1$ or
$\kappa(D\vdual[p])<1$. Here $D\vdual$ is the Cartier
dual of $D$. We say the condition $\Cond$ holds for an
$m$-truncated Barsotti--Tate group $B$ over a field $K$ of
characteristic $p$, if the condition $\Cond$ holds for the
extension of $B$ to an algebraic closure of $K$.
\end{proclamation}

\begin{proclamation}{Lemma}\label{CompLem}
The following four properties hold:

\medskip
\begin{definitionlist}
\item Let $\gamma\in\Gamma$. Then for each positive integer $t$ we have the following inequality $n_t(\gamma)\le 1+\sum_{u=1}^{t-1} n_u(\gamma).$ Therefore we have $n_t(\gamma)\le 2^{t-1}$.

\item
If $p\ge 3$ and $\gamma\in\Gamma$, then $\kappa(\gamma)<1$.

\item
If  $p\ge 5$ and $\gamma\in\Gamma\cup\Delta$, then
$\kappa(\gamma)<1$. Thus if $p\ge 5$, then $\kappa(\pi)<1$ and
therefore condition $\Cond$ holds for $D$.

\item
We assume that $p=3$. Then $\kappa(\pi)<{4\over 3}$. If moreover $\min\{c,d\}\le 6$, then condition $\Cond$ holds for $D$.
\end{definitionlist}
\end{proclamation}

\begin{proof}
To prove (a) we write $\gamma=(i_1,\ldots,i_s)$ and we can assume that
$n_t(\gamma)>0$; thus we have $s\ge 3$ even if $\gamma\in\Gamma$. We consider the identity
$$\bar e_{\pi^t(i_1),\pi^t(i_{s-1})}=\prod_{\ell=1}^{s-2} \bar
e_{\pi^t(i_{\ell}),\pi^t(i_{\ell+1})}\leqno (1)$$ between elements
of $\End(M/pM)$. The right hand side of (1) contains at least
$n_t(\gamma)$ elements of $\nfr_+$, it contains $\sum_{u> t}
n_u(\gamma)$ elements of $\nfr_{0,0}$, and it contains
$n_1(\gamma)+n_2(\gamma)+\cdots+n_{t-1}(\gamma)$ elements that
belong to $\Jscr_+$, $\Jscr_0$, or $\Jscr_-$. Let
$\bar\mu:\Bbb G_m\to\pmb{GL}_{M/pM}$ be the cocharacter that fixes
$\bar F^0$ and that acts via the inverse of the identical
character of $\Bbb G_m$ on $\bar F^1$. If we have
$n_t(\gamma)>1+n_1(\gamma)+n_2(\gamma)+\cdots+n_{t-1}(\gamma)$,
then $\Bbb G_m$ acts via $\bar\mu$ on the right hand side of (1)
via at least the second power of the inverse of the identical
character of $\Bbb G_m$; therefore the right hand side of (1) must
be $0$ and this contradicts the fact that the left hand side of
(1) is non-zero. Thus we have $n_t(\gamma)\le
1+n_1(\gamma)+n_2(\gamma)+\cdots+n_{t-1}(\gamma)$. By induction on  
$t\in\Bbb N^*$ we
get that $n_t(\gamma)\le 2^{t-1}=1+1+2+\cdots +2^{t-2}$. This
proves (a).

We prove (b). Due to (a) and the inequality $p\ge 3$ we have
$$\kappa(\gamma)=\sum_{t=1}^{\infty} n_t(\gamma)p^{-t}< \sum_{t=1}^{\infty}
2^{t-1}p^{-t}={1\over {p-2}}\le 1.\leqno (2)$$
Thus (b) holds.

Due to (b), to prove (c) we can assume that
$\gamma=(i_1,\ldots,i_s)\in \Delta$. We have
$$\kappa(\gamma)={1\over
{p^{\nu(i_{s-1},i_s)}}}+\kappa((i_1,\ldots,i_{s-1}))\le {1\over
p}+\kappa((i_1,\ldots,i_{s-1})).\leqno (3)$$ Due to the last
inequality and the fact that $p\ge 5$, from relations (2) and (3) we get
that $\kappa(\gamma)<{1\over {p-2}}+{1\over p}<1$. Thus (c) holds.

The first part of (d) follows from (b) and (3). To check the last part of (d) we can assume that $d\le 6$ (otherwise we can replace
$D$ by $D\vdual$). From Lemma \ref{NilLem} (c), (by taking $t=0$ in Formula (1)) we get
that for each element $\gamma=(i_1,\ldots,i_s)\in\Gamma$ we have
$s\le d+1\le 7$ and thus $\sum_{t=1}^{\infty}
n_t(\gamma)\le s-2\le 5$. From this and (a) we get that
$\kappa(\gamma)\le {1\over 3}+{2\over 9}+{2\over 27}={{17}\over
{27}}$. As ${{17}\over {27}}+{1\over 3}={{26}\over {27}}<1$, from the last sentence
and (3), we get that we have $\kappa(\tilde\gamma)<1$ for all
$\tilde\gamma\in\Delta$. Thus $\kappa(\pi)<1$ and therefore (d)
holds.
\end{proof}

\begin{remark}{Example}\label{Ex1}
We assume that $p\in \{2,3\}$ and that there exists an integer $a\geq 2$
such that we have a ring monomorphism $\Bbb
F_{p^a}\hookrightarrow \End(D[p])$ with the property that $\Bbb
F_{p^a}$ acts on the tangent space of $D[p]$ via scalar
endomorphisms. We will check that $\kappa(D)<1$.

To the product decomposition $\Bbb F_{p^a} \otimes_{\Bbb
F_p} k=k^a$ corresponds a direct sum decomposition
$M/pM=\oplus_{u=1}^a \bar M_u$. It is easy to see that we can
assume that the last sum decomposition is the reduction modulo $p$
of a direct sum decomposition $M=\oplus_{u=1}^a M_u$ with the
property that $\sigma_{\pi}(M_u)=M_{u+1}$, where $M_{a+1}:=M_1$.
The direct sum decomposition  of $W(k)$-modules  
$\End(M)=\oplus_{u=1}^a\oplus_{\tilde u=1}^a \Hom(M_u,M_{\tilde u})$  
allows us to
view each $\Hom(M_u,M_{\tilde u})$ as a $W(k)$-submodule of
$\End(M)$.

\noindent As $\Bbb F_{p^a}$ acts on the tangent space of $D[p]$
via scalar endomorphisms, we can choose the indexing of $\bar
M_u$'s such that we have $\bar F^1\subseteq \bar M_1$. Thus
we can assume that $F^1\subseteq M_1$. We can also assume
that certain subsets of $\{e_1,\ldots,e_r\}$ are $W(k)$-bases for
the $M_u$'s. From the last two sentences, we get that:

\medskip
{\bf (*)}  if $(i,j)\in\Jscr_+$, then $e_{i,j}\in  
\Hom(M,F^1)\subseteq\Hom(M,M_1)$.

\medskip
\noindent Let $\gamma=(i_1,\ldots,i_s)\in\Gamma\cup\Delta$. To
check that $\kappa(D)<1$, it suffices to show that
$\kappa(\gamma)<1$. Based on (3), to check this it suffices to
show that if $\gamma\in\Gamma$, then we have
$\kappa(\gamma)<1-{1\over p}$. We can assume that $s\ge 3$. As
$\gamma\in\Gamma$, we have $(i_{s-1},i_s)\in\Jscr_{+,2}$ and
$(i_1,i_2),\ldots,(i_{s-2},i_{s-1})\in\Jscr_{0,0}$. From this and
(*) we easily get that we have
$e_{i_1,i_2},\ldots,e_{i_{s-2},i_{s-1}}\in\End(M_1)$. As
$\sigma_{\pi}(\End(M_u))=\End(M_{u+1})$, we get that each positive
integer $i$ such that $e_{\pi^i(i_1),\pi^i(i_2)}\in\End(M_1)$ is a
multiple of $a$. From this and (*), we get that $\nu(i_1,i_2)\subseteq a\Bbb N^*$. Similarly we argue that $\{\nu(i_2,i_3),\ldots,\nu(i_{s-2},i_{s-1})\}\subseteq a\Bbb N^*$. Therefore for each $t\in\Bbb N^*\setminus a\Bbb N^*$ we have
$n_t(\gamma)=0$. As in the proof of Lemma \ref{CompLem} (a) we
argue that for each $t\in\Bbb N^*$ we have $n_{ta}(\gamma)\le
2^{t-1}$. Therefore we have $\kappa(\gamma)<\sum_{t=1}^{\infty}
2^{t-1}p^{-at}={1\over {p^a-2}}$. As $a\ge 2$, we get that
$\kappa(\gamma)< 1-{1\over p}$. Thus $\kappa(D)<1$.
\end{remark}

\begin{remark}{Remark}\label{Badcase}
If $p=2$ (resp. $p=3$), there exist plenty of examples in which there  
exist elements $\gamma\in\Gamma_1$ (resp. $\gamma\in\Delta_1$) such  
that $\kappa(\gamma)>1$ and therefore also $\kappa(\pi)>1$ (for $p=2$  
see Example \ref{Ex3} below).
\end{remark}
\end{segment}

\begin{segment}{Computing $\Cscr_1^0$ with explicit equations}{CompC10}
Let $(h_1[1],h_2[1],h_3[1])\in\Hscr_1(k)$. Let
$h_{12}[1]:=h_1[1]h_2[1]\in\Wscr_{+0}(k)$ and
$h_{23}[1]:=h_2[1]h_3[1]\in\Wscr_{0-}(k)$. We have $(h_1[1],h_2[1],h_3[1])\in\Cscr_1(k)$ if and only if $h_{12}[1]=\sigma_{\phi}(h_{23}[1])$, cf. the very definition of the action $\Bbb T_1$. Writing $h_{12}=1_M[1]+\sum_{(i,j)\in\Jscr_+\cup\Jscr_0} x_{i,j}\bar e_{i,j}$ and $h_{23}=1_M[1]+\sum_{(i,j)\in\Jscr_0\cup\Jscr_-} x_{i,j}\bar e_{i,j}$ with $x_{i,j}$ as independent variables, the equation $h_{12}[1]=\sigma_{\phi}(h_{23}[1])$ gets translated into a system of equations that are of the form $a_{\pi(i),\pi(j)}x_{\pi(i),\pi(j)}=b_{i,j}x^p_{i,j}$ with $a_{\pi(i),\pi(j)},b_{i,j}\in\{0,1\}$ and that are indexed by $(i,j)\in J^2$. More precisely, we have $a_{i,j}=1$ (resp. $b_{i,j}=1$) if and only if $(i,j)\in \Jscr_+\cup\Jscr_0$ (resp. $(i,j)\in \Jscr_0\cup\Jscr_-$). Based on this, by tracing those variables $x_{i,j}$ that can take independently an infinite number of values in $k$, one easily gets that we have
$(h_1[1],h_2[1],h_3[1])\in\Cscr_1^0(k)$ if and only if the following three identities hold (to be compared with \cite{Va2}, Subsection 2.3, Formulas $(4a)$ to
$(4c)$): 
$$h_{12}[1]=1_M[1]+\sum_{(i,j)\in\Jscr_{-,1}} \sum_{\ell=1}^{\nu(i,j)}
x_{i,j}^{p^\ell}\bar e_{\pi^\ell(i),\pi^\ell(j)},\leqno (4)$$
$$h_2[1]=1_M[1]+\sum_{(i,j)\in\Jscr_{-,1}}\sum_{\ell=1}^{\nu(i,j)-1}
x_{i,j}^{p^\ell}\bar e_{\pi^\ell(i),\pi^\ell(j)},\leqno (5)$$
$$h_{23}[1]=1_M[1]+\sum_{(i,j)\in\Jscr_{-,1}} \sum_{\ell=0}^{\nu(i,j)-1}
x_{i,j}^{p^\ell}\bar e_{\pi^\ell(i),\pi^\ell(j)},\leqno (6)$$
where the variables $x_{i,j}$ with $(i,j)\in\Jscr_{-,1}$ can take
independently all values in $k$ such that $h_2[1]\in\Wscr_0(k)$.

Note that formulas $(4),(5)$, and $(6)$ differ only by
summation limits. If $\nfr_{0,0}$ is as in Lemma \ref{NilLem} (b),
then we have $h_2[1]\in 1_M[1]+\nfr_{0,0}$. From this and the fact
that $\nfr_{0,0}$ is a nilpotent algebra, we get that:

\medskip
{\bf (*)} {\it the element $h_2[1]$ has an inverse in  
$1_M[1]+\nfr_{0,0}$ which is a
polynomial in $h_2[1]$ with coefficients in $k$ and therefore we always have
$h_2[1]\in\Wscr_0(k)$.}

\medskip\noindent
Thus the variables $x_{i,j}$ with $(i,j)\in\Jscr_{-,1}$ can take
independently all values in $k$. Based on Formulas $(4)$ to
$(6)$ we get that
$$\Lie(\Cscr_1^0)=\bigoplus_{(i,j)\in\Jscr_{-,1}} k\bar
e_{i,j}=\nfr_{-,1}\subseteq\Lie(\Wscr_{-,k}).\leqno(7)$$
As $\nfr_{-,1}$ does not contain non-zero semisimple elements, from (7) we get that $\dbG_m$ is not a subgroup of $\Cscr_1^0$ and thus $\scrC_1^0$ is unipotent.

Let $\Vscr_{-1}$, $\Vscr_1$, and $\Vscr_2$ be the vector group
schemes over $k$ whose Lie algebras are $\nfr_{-,1}$, $\nfr_{+,1}$
and $\nfr_{+,2}$ (respectively). For a commutative $k$-algebra $R$
we have $\Vscr_{-1}(R)=\nfr_{-,1}\otimes_k R$,
$\Vscr_1(R)=\nfr_{+,1}\otimes_k R$, and
$\Vscr_2(R)=\nfr_{+,2}\otimes_k R$. We get that the morphism of smooth $k$-schemes 
$$\Vscr_{-1}\to \Cscr_1^0$$
that takes the element $\sum_{(i,j)\in\Jscr_{-,1}} x_{i,j}\bar e_{i,j}\in
\nfr_{-,1}=\Lie(\Cscr_1^0)=\Lie(\Vscr_{-1})=\Vscr_{-1}(k)$ to the
element
$$(h_{12}[1](h_2[1])^{-1},h_2[1],(h_2[1])^{-1}h_{23}[1])\in\Cscr^0_1(k)$$
obtained naturally from formulas $(4)$ to $(6)$, is an
isomorphism. For $s\in\{1,2\}$ we consider the closed embedding
monomorphism
$$\varepsilon_s:\Vscr_s\hookrightarrow\Wscr_{+,k}$$
that takes $v_s\in \Vscr_s(R)$ to $1_{M/pM\otimes_k R} + v_s\in  
\Wscr_{+,k}(R)$.
\end{segment}

\begin{segment}{The key morphism $\zeta$}{KeyZeta}
Let $\Iscr_1:=\Hscr_1/\Wscr_{-,k}$; it is an affine, smooth group
scheme over $k$ isomorphic to $\Wscr_{+0,k}$ and therefore such
that we have a short exact sequence
$$1\to \Wscr_{+,k}\to \Iscr_1\to\Wscr_{0,k}\to 1.$$
As schemes over $k$ we can identify naturally
$\Iscr_1=\Wscr_{+,k}\times_k \Wscr_{0,k}$. Thus we will also
identify (as sets) $\Iscr_1(k)=\Wscr_+(k)\times \Wscr_0(k)$. Let
$\epsilon:\Hscr_1\twoheadrightarrow \Iscr_1$ and
$\theta:\Iscr_1\twoheadrightarrow \Wscr_{0,k}$ be the natural
epimorphisms and let $\iota:=\theta\circ\epsilon:\Hscr_1\twoheadrightarrow
\Wscr_{0,k}$ be their composite. The epimorphism $\epsilon$ gives birth via restriction
to a finite homomorphism $\Cscr_{1}^0\to \Iscr_{1}$ whose kernel
is finite and connected. 

The group $\tilde\Cscr_1^0 :=
\text{Im}(\Cscr_{1}^0\to\Iscr_{1})$ is isomorphic to $\Cscr_1^0$.
More precisely, using the isomorphism $\Vscr_{-1}\to\scrC_1^0$ one easily gets that the epimorphism
$\Cscr_1^0\twoheadrightarrow\tilde\Cscr_1^0$ can be identified
with the Frobenius endomorphism of $\Cscr_1^0$. We have
$$\Lie(\tilde\Cscr_1^0)=\bigoplus_{(i,j)\in\Jscr_{0,0}\cup\Jscr_{+,1},\eta(i,j)=1}\;
k\bar e_{i,j}.$$

\begin{proclamation}{Definition}
The morphism $\zeta:\tilde\Cscr_{1}^0\times_k \Vscr_2\times_{k}
\Wscr_{0,k}\to \Iscr_1$ of $k$-schemes is defined by the following
rule: if $R$ is a commutative $k$-algebra, then the element
$(\tilde h,y,z)\in
\tilde\Cscr_1^0(R)\times\Vscr_2(R)\times\Wscr_{0,k}(R)$ is mapped
to the product element $\tilde h\cdot (\varepsilon_2(y),1)\cdot
\epsilon(1,z,1)\in \Iscr_1(R).$
\end{proclamation}

\noindent The key step for proving the main result
Theorem~\ref{MainThm} is to show that $\zeta$ is finite and flat.
For this, we first prove the following basic fact.

\begin{proclamation}{Lemma}\label{GrobLem}
Let $R$ be a commutative $\Bbb F_p$-algebra. Let $n\in\Bbb N^*$. Let
$d_1,\ldots,d_n$ be positive integers. For each
$\ell\in\{1,\ldots,n\}$ we consider a polynomial $Q_{\ell}\in
R[x_1,\ldots,x_n]$. We assume that there exist positive rational numbers
$\mu_1,\ldots,\mu_n$ such that for each $\ell\in\{1,\ldots,n\}$
and for every monomial term $\beta\prod_{i=1}^n x_i^{v_i}$ of
$Q_{\ell}$ with $\beta\in R\setminus\{0\}$ we have an inequality  
$\mu_{\ell} d_{\ell} >\sum_{i=1}^n
\mu_iv_i$ (thus the degree
of $Q_{\ell}$ in each variable $x_i$ is at most equal to $d_i-1$). Let
$$\Rscr:=R[x_1, \ldots,x_n]/(x_1^{d_1}-Q_1,\ldots,x_n^{d_n}-Q_n).$$
Then $\Rscr$ is a free $R$-module of rank $\prod_{\ell=1}^n d_{\ell}$.
Therefore the natural morphism $\Spec \Rscr\to\Spec R$ of schemes is a
finite, flat cover of degree $\prod_{\ell=1}^n d_{\ell}$.
\end{proclamation}

\begin{proof}
Multiplying all $\mu_1,\ldots,\mu_n$ by a positive integer, we can  
assume that we have $\mu_1,\ldots,\mu_n\in\Bbb N^*$. Let $\Bbb M:=\Bbb N^d$. Let $\tau:\Bbb
M\hookrightarrow\Bbb M$ be the monomorphism of additive monoids
that takes a sequence $(v_1,\ldots,v_n)\in\Bbb M$ to
$(\mu_1v_1,\ldots,\mu_nv_n)\in\Bbb M$. Let $<$ be the
degree-lexicographic ordering on $\Bbb M$. Let $\prec$ be the well
ordering on $\Bbb M$ such that for $a,b\in\Bbb M$ we have $a\prec
b$ if and only if $\tau(a)<\tau(b)$. For each $P\in
R[x_1, \ldots,x_n]$, we define its weight $\omega(P)$ to be the
maximal element $(v_1,\ldots,v_n)$ of $\Bbb M$ with respect to the
ordering $\prec$ for which the monomial $\prod_{i=1}^n x_i^{v_i}$
appears in $P$ with a non-zero coefficient. From hypotheses we get
that for each $\ell\in\{1,\ldots,n\}$ we have
$\omega(x_{\ell}^{d_{\ell}})>\omega(Q_{\ell})$. Using this, it is an easy
exercise in the theory of Gr\"obner bases (over the base $R$) to
check that the image of the set $\{\prod_{i=1}^n x_i^{v_i}|0\le
v_i<d_i \;\;\forall\;\; i\in\{1,\ldots,n\}\}$ in $\Rscr$ is an
$R$-basis for $\Rscr$. From this the Lemma follows.
\end{proof}

\begin{proclamation}{Theorem}\label{FinThm}
We assume that we have an inequality $\kappa(D[p])<1$. Then the  
morphism $\zeta$ is
a finite, flat morphism of degree $p^{|\Jscr_{0,0}|}$.
\end{proclamation}

\begin{proof} Each element of $\Wscr_+(k)$ can be written uniquely as
$\varepsilon_1(a_1)\varepsilon_2(a_2)$ with $a_1\in\nfr_{+,1}$ and
$a_2\in\nfr_{+,2}$. Thus each element of $\Iscr_1(k)$ can be
written uniquely as a pair $(\varepsilon_1(a_1)\varepsilon_2(a_2),a_0)$
with $a_0\in\Wscr_0(k)$. Let $(\tilde h,y,z)\in \tilde\Cscr_1^0(k) \times
\Vscr_2(k)\times\Wscr_{0,k}(k)$ be
an arbitrary element. We look at the equation
$$\zeta((\tilde h,y,z))=(\varepsilon_1(a_1)\varepsilon_2(a_2),a_0).\leqno (8)$$
By applying $\theta$ to $(8)$ we get that
$$\iota(\tilde h)z=a_0.\leqno (9)$$
We write $\tilde h=(h_1(x),h_2(x))$ with
$x\in\Lie(\tilde\Cscr_1^0)$. We have $h_2(x)z=a_0$, cf. (9).

\noindent
We write $y=\sum_{(i,j)\in\Jscr_{+,2}} y_{i,j}\bar e_{i,j}$ and
$x=\sum_{(i,j)\in\Jscr_{0,0}\cup\Jscr_{+,1},\eta(i,j)=1}
x_{i,j}\bar e_{i,j}$, where $y_{i,j}$'s and $x_{i,j}$'s are
variables. Based on formulas $(4)$ and $(5)$, we have
$$h_2(x)=1_M[1]+\sum_{(i,j)\in\Jscr_{0,0},\eta(i,j)=1}\sum_{l=0}^{\nu(i,j)-1}
x_{i,j}^{p^l}\bar e_{\pi^l(i),\pi^l(j)}$$
and
$$h_{12}(x):=h_1(x)h_2(x)=1_M[1]+\sum_{(i,j)\in\Jscr_{0,0}\cup\Jscr_{+,1},\eta(i,j)=1}\sum_{l=0}^{\nu(i,j)}
x_{i,j}^{p^l}\bar e_{\pi^l(i),\pi^l(j)}.$$
By applying $(8)$ and $(9)$, by using the identification
$\Iscr_1(k)=\Wscr_+(k)\times\Wscr_0(k)$, and by denoting with $1$ identity elements, we get that

$$(h_1(x),h_2(x))\cdot (1+y,1)\cdot (1,z) = (1+a_1+a_2,1)\cdot (1,h_2(x))\cdot
(1,z).$$

\noindent Thus we get the equation $(h_1(x),h_2(x))\cdot
(1+y,1)=(1+a_1+a_2,1)\cdot (1,h_2(x))$ between elements of
$\Iscr_1(k)$ and therefore the equation
$$h_{12}(x)(1_M[1]+y)=(1_M[1]+a_1+a_2)h_2(x)\leqno (10)$$
between elements of $\pmb{\text{GL}}_M(k)$.  As $\nfr_+^2=0$, we
have $h_{12}(x)y=h_2(x)y$. Thus Equation $(10)$ is equivalent to $h_{12}(x)+h_2(x)y=h_2(x)+(a_1+a_2)h_2(x)$. Let
$$l_1(x):=h_{12}(x)-h_2(x)=\sum_{(i,j)\in\Jscr_{0,0}\cup\Jscr_{+,1},\eta(i,j)=1}
x_{i,j}^{p^{\nu(i,j)}}\bar  
e_{\pi^{\nu(i,j)}(i),\pi^{\nu(i,j)}(j)}\in\nfr_{+,1}$$
and
$$l_2(x):=h_2(x)-1_M[1]=\sum_{(i,j)\in\Jscr_{0,0},\eta(i,j)=1}\sum_{l=0}^{\nu(i,j)-1}
x_{i,j}^{p^l}\bar
e_{\pi^l(i),\pi^l(j)}=:\sum_{(i^\prime,j^\prime)\in\Jscr_{0,0}}
q_{i^\prime,j^\prime}(x)\bar e_{i^\prime,j^\prime}.$$
We get that
Equation $(10)$ can be rewritten as
$$l_1(x)+h_2(x)y=(a_1+a_2)h_2(x).\leqno (11)$$
For $s\in\{1,2\}$ let
$$\pi_s:\Lie(\Wscr_{+,k})\twoheadrightarrow\Lie(\Vscr_s)$$
be the projection of $\Lie(\Wscr_{+,k})$ on $\Lie(\Vscr_s)$ along
$\Lie(\Vscr_{3-s})$. Multiplying $(11)$ from the left by $h_2(x)^{-1}$  
and using
that $h_2(x)^{-1}l_1(x)=\sum_{s=0}^{\infty} (-1)^sl_2(x)^sl_1(x)\in  
\nfr_{+,1}$ (cf. property 2.5 (*) and the Lemma \ref{NilLem} (c)), we get that
$$y=\pi_2\bigl(h_2(x)^{-1}a_2h_2(x)\bigr).\leqno (12)$$
Due to $(12)$, solving Equation $(11)$ is the same thing as solving the
equation $h_2(x)^{-1}l_1(x)=\pi_1\bigl(h_2(x)^{-1}(a_1+a_2)h_2(x)\bigr)$ and
therefore the equation
$$l_1(x)=[1_M[1]+l_2(x)]\pi_1\bigl(\sum_{s=0}^{\infty}
(-1)^sl_2(x)^s(a_1+a_2)[1_M[1]+l_2(x)]\bigr)=$$
$$a_1+a_1l_2(x)+[1_M[1]+l_2(x)]\pi_1\bigl(\sum_{s=0}^{\infty}
(-1)^sl_2(x)^s(a_2)[1_M[1]+l_2(x)]\bigr)
.\leqno (13)$$
The last identity is implied by Lemma \ref{NilLem} (c). 

\noindent We write $a_1+a_2=\sum_{(i,j)\in\Jscr_+} a_{i,j}\bar
e_{i,j}$, with the $a_{i,j}$'s thought as variables. We consider
the polynomial $k$-algebra $R:=k[a_{i,j}|(i,j)\in\Jscr_+]$. Due to
the above formulas for $l_1(x)$ and $l_2(x)$, the system of
equations defined by Equation $(13)$ has the form
$$x_{i,j}^{p^{\nu(i,j)}}=Q_{i,j}\;\;\;\text{for}\;\;\;
(i,j)\in\Jscr_{0,0}\cup\Jscr_{+,1}\;\;\;\text{with}\;\;\;\eta(i,j)=1,\leqno
(14)$$ where each $Q_{i,j}$ is a polynomial in the variables
$x_{i,j}$'s with coefficients in $R$. This system defines a
morphism $\Spec \Rscr\to\Spec R$ of affine $k$-schemes. For $(i,j)\in\Jscr_{0,0}\cup\Jscr_{+,1}$ with $\eta(i,j)=1$, by identifying the coefficients of $\bar e_{\pi^{\nu(i,j)}(i),\pi^{\nu(i,j)}(j)}$ in the two sides of Equation (13), we get the formula:
$$Q_{i,j}=\sum_{(i^\prime_1,i^\prime_2,\ldots,i^\prime_{s-1},i^\prime_s)\in\Gamma_1,\;i^\prime_1=\pi^{\nu(i,j)}(i),\;i^\prime_s=\pi^{\nu(i,j)}(j)}
(-1)^{s-2}a_{i^\prime_{s-1},i^\prime_s}q_{i^\prime_1,i^\prime_2}(x)\cdots  
q_{i^\prime_{s-2},i^\prime_{s-1}}(x)+$$
$$\sum_{(i^\prime_1,i^\prime_2,\ldots,i^\prime_{s-1},i^\prime_s)\in\Delta_1,\;i^\prime_1=\pi^{\nu(i,j)}(i),\;i^\prime_s=\pi^{\nu(i,j)}(j)}
(-1)^{s-3}a_{i^\prime_{s-2},i^\prime_{s-1}}q_{i^\prime_1,i^\prime_2}(x)\cdots
q_{i^\prime_{s-3},i^\prime_{s-2}}(x)q_{i^\prime_{s-1},i^\prime_s}(x)+$$
$$a_{\pi^{\nu(i,j)}(i),\pi^{\nu(i,j)}(j)}+\sum_{i'\in\{1,\ldots,r\},(\pi^{\nu(i,j)}(i),i')\in\scrJ_{+,1},(i',\pi^{\nu(i,j)}(j))\in\scrJ_{0,0}}a_{\pi^{\nu(i,j)}(i),i'}q_{i',\pi^{\nu(i,j)}(j)}(x).\leqno  
(15)$$
In Formula (15), the first two lines  record (resp. the third line records) the contribution of $a_2$ (resp. of $a_1$) to $Q_{i,j}$. 

\noindent Let $n$ be the number of elements of $\Jscr_{-,1}$. We
write
$$\{(i,j)\in\Jscr_{0,0}\cup\Jscr_{+,1}|\eta(i,j)=1\}=\{(i_1,j_1),\ldots,(i_n,j_n)\}.$$
  We can assume that $\kappa(\pi)<1$, cf. the very definition of
$\kappa(D[p])$. Let
$$\nu:=\max\{\nu(i_{\ell},j_{\ell})|\ell\in\{1,\ldots,n\}\}.$$
For each
$\ell\in\{1,\ldots,n\}$ we define a positive integer (thought as a  
weight) $\mu_{\ell}:=p^{\nu-\nu(i_{\ell},j_{\ell})}$. We
have $\mu_{\ell}p^{\nu(i_{\ell},j_{\ell})}=p^{\nu}$.

\noindent We consider a monomial term:
$$\beta\prod_{\ell=1}^n  x_{i_{\ell},j_{\ell}}^{t_{\ell}}=(-1)^{s-2}a_{i^\prime_{s-1},i^\prime_s}q_{i^\prime_1,i^\prime_2}(x)\cdots
 q_{i^\prime_{s-2},i^\prime_{s-1}}(x)$$
  with  
$\gamma:=(i^\prime_1,i^\prime_2,\ldots,i^\prime_{s-1},i^\prime_s)\in\Gamma_1$  
(resp.
$$\beta\prod_{\ell=1}^n  x_{i_{\ell},j_{\ell}}^{t_{\ell}}=(-1)^{s-3}a_{i^\prime_{s-2},i^\prime_{s-1}}q_{i^\prime_1,i^\prime_2}(x)\cdots
 q_{i^\prime_{s-3},i^\prime_{s-2}}(x)q_{i^\prime_{s-1},i^\prime_s}(x)$$
with   
$\gamma:=(i^\prime_1,i^\prime_2,\ldots,i^\prime_{s-1},i^\prime_s)\in\Delta_1$)  
that
contributes via the Equation $(15)$ to some $Q_{i_a,j_a}$. If
$u\in\{1,\ldots,s-1\}$ is such that $(i^\prime_u,i^\prime_{u+1})\in  
\Jscr_{0,0}$, then for  
$(i_{\ell},j_{\ell}):=\pi^{-\eta(i^\prime_u,i^\prime_{u+1})+1}(i^\prime_u,i^\prime_{u+1})\in\{(i,j)\in
\Jscr_{0,0}\cup\Jscr_{+,1}|\eta(i,j)=1\}$ with $\ell\in\{1,\ldots,n\}$ we get that the contribution of $x_{i_{\ell},j_{\ell}}$ to $\beta\prod_{\ell=1}^n
x_{i_{\ell},j_{\ell}}^{t_{\ell}}$ that corresponds to the segment  
$(i^\prime_u,i^\prime_{u+1})$ of $\gamma$ is  
$x_{i_{\ell},j_{\ell}}^{p^{\nu(i_{\ell},j_{\ell})-\nu(i^\prime_u,i^\prime_{u+1})}}$.
Thus the contribution of  
$x_{i_{\ell},j_{\ell}}$
to the sum $\sum_{\ell=1}^n t_{\ell}\mu_{\ell}$ that corresponds to  
the segment $(i^\prime_u,i^\prime_{u+1})$  is precisely
$p^{\nu-\nu(i^\prime_u,i^\prime_{u+1})}$. For $t\in \Bbb N^*$, let
$n_t(\gamma)$ be as in Subsection \eqref{COND} (i.e., the number
of those $u\in\{1,\ldots,s-1\}$ such that we have  
$(i^\prime_u,i^\prime_{u+1})\in \Jscr_{0,0}$ and  
$\nu(i^\prime_u,i^\prime_{u+1})=t$). From the last two sentences, we  
get that
$$\sum_{\ell=1}^n t_{\ell}\mu_{\ell}=p^{\nu}\sum_{t=1}^{\infty}  
n_t(\gamma)p^{-t}=p^{\nu}\kappa(\gamma)\le
p^{\nu}\kappa(\pi)<p^{\nu}\leqno (16)$$
(the first inequality, cf. the very definition of $\kappa(\pi)$). 

\noindent
If $\beta\prod_{\ell=1}^n
x_{i_{\ell},j_{\ell}}^{t_{\ell}}$ is $a_{\pi^{\nu(i,j)}(i),\pi^{\nu(i,j)}(j)}$  (resp. is $a_{\pi^{\nu(i,j)}(i),i'}q_{i',\pi^{\nu(i,j)}(j)}(x)$ with $\break(\pi^{\nu(i,j)}(i),i')\in\scrJ_{+,1}$ and $(i',\pi^{\nu(i,j)}(j))\in\scrJ_{0,0}$), then the sum $\sum_{\ell=1}^n t_{\ell}\mu_{\ell}$ is $0$ (resp. is $p^{\nu-\nu(i',\pi^{\nu(i,j)}(j))}\le p^{\nu-1}$) and thus it is less than $p^{\nu}$. 

\noindent
From the last two paragraphs we get that for every monomial term  
$\beta\prod_{\ell=1}^n
x_{i_{\ell},j_{\ell}}^{t_{\ell}}$ that shows up with a non-zero  
coefficient in some $Q_{i_a,j_a}$ with $a\in\{1,\ldots,n\}$, we have $\sum_{\ell=1}^n t_{\ell}\mu_{\ell}<p^{\nu}$ and thus (cf. (14)) the
hypotheses of Lemma \ref{GrobLem} are satisfied. Thus from Lemma  
\ref{GrobLem} we get
that the morphism $\Spec \Rscr\to\Spec R$ is a finite, flat cover of degree
$\prod_{\ell=1}^n  
p^{\nu(i_{\ell},j_{\ell})}=p^{\sum_{(i,j)\in\Jscr_{-,1}}  
\nu(i,j)-1}=p^{|\Jscr_{0,0}|}$.

All the above part can be redone using arbitrary $\Kscr$-valued
points instead of $k$-valued points. Taking $\Kscr$ to be the
$k$-algebra of global functions of the affine $k$-scheme
$\Iscr_1$, from Equation $(9)$ and the fact that the morphism
$\Spec \Rscr\to\Spec R$ is a finite, flat cover of degree
$p^{|\Jscr_{0,0}|}$, we get that $\zeta$ is a finite, flat cover
of degree $p^{|\Jscr_{0,0}|}$.
\end{proof}
\end{segment}

\begin{segment}{Examples}{Examp}
\begin{remark}{Example}\label{Ex2}
Suppose that $D$ is minimal i.e., we have $n_D=1$. As $n_D=1$, we
can assume that $g_{\pi}=1_M$ and thus that
$\sigma_{\phi}=\sigma_{\pi}$. As $\sigma_{\phi}=\sigma_{\pi}$, we have  
direct sum
decompositions:
$$(\End(M),\phi)=\oplus_{\alpha\in \Bbb Q \cap [-1,1]}
(W_{\alpha},\phi)=(S_+,\phi)\oplus (S_0,\phi)\oplus (S_-,\phi)$$
such that all Newton polygon slopes of $(W_{\alpha},\phi)$ are
$\alpha$, all Newton polygon slopes of $(S_+,\phi)$ are positive,
all  Newton polygon slopes of $(S_0,\phi)$ are $0$, and all Newton
polygon slopes of $(S_-,\phi)$ are negative. 

Referring to Equation
$(11)$,  in this paragraph we will check that $y,a_2\in S_+/pS_+$. To check this we
can assume that either $D$ is indecomposable or a direct sum of
two indecomposable minimal $p$-divisible groups over $k$ that have
distinct Newton polygon slopes. In the first case we have
$S_+=S_-=0$ and from \cite{Va2}, Example 2.3.1 we get that
$\Jscr_{-,2}=\Jscr_-\setminus\Jscr_{-,1}=\emptyset$ and therefore
that $\Jscr_{+,1}=\Jscr_+$ and $\Jscr_{+,2}=\emptyset$; thus
$y=a_2=0\in S_+/pS_+$. In the second case, we write $D = D_1
\times D_2$ such that both $D_1$ and $D_2$ are indecomposable
minimal $p$-divisible groups over $k$. For $s\in\{1,2\}$, let
$c_s$ and $d_s$ be the codimension and the dimension of $D_s$. We
can assume that the Newton polygon slope $\frac{d_1}{c_1+d_1}$ of
$D_1$ is less than the Newton polygon slope $\frac{d_2}{c_2+d_2}$
of $D_2$. We have $c = c_1 + c_2, d = d_1 + d_2$, and $c_2d_1 <
c_1d_2$. Let $(M,\phi)=(M_1,\phi)\oplus (M_2,\phi)$ be the direct
sum decomposition such that the Dieudonn\'e module of $D_s$ is
$(M_s,\phi)$. We have $S_+=\Hom(M_1,M_2)$,
$S_0=\End(M_1)\oplus\End(M_2)$, and $S_-=\Hom(M_2,M_1)$. From
\cite{Va2}, Example 2.3.2, we get that if $(i,j)\in\Jscr_{-}$ is
such that $e_{i,j}\in \Hom(M_1,M_2)$, then we have
$(i,j)\in\Jscr_{-,1}$. The $W(k)$-linear map that takes $e_{i,j}$
to $e_{j,i}$ induces naturally an isomorphism between
$(\Hom(M_2,M_1),\phi)$ and the dual of $(\Hom(M_1,M_2,\phi))$ (in
the natural sense that involves no Tate twists). From the previous
two sentences we get via duality that if $(i,j)\in\Jscr_+$ is such
that $e_{i,j}\in \Hom(M_2,M_1)$, then we have
$(i,j)\in\Jscr_{+,1}$. From this and the first case, we get that
if $(i,j)\in\Jscr_{+,2}$, then we have $e_{i,j}\in\Hom(M_1,M_2)=S_+$.
This implies that $y,a_2\in S_+/pS_+$.

In this paragraph we assume that $a_1=0$. As $y,a_2\in
S_+/pS_+$ and $a_1=0$, from Equation $(11)$ we get by an increasing induction
on the Newton polygon slope $\alpha\in \Bbb Q \cap [-1,1]$ that the
component of $x$ in $W_{\alpha}/pW_{\alpha}$ is uniquely
determined. Thus the system $(11)$ has for $a_1=0$ a unique solution.
More precisely, for $a_1=0$ the system $(13)$ is of the form:
$$x_1^{v_1}=a_1,\;\;\;\;x_2^{v_2}=Q_2(x_1),\;\;\;\;\cdots,\;\;\;\;x_n^{v_n}=Q_n(x_1,\ldots,x_{n-1}),$$
where $\{x_1,\ldots,x_n\}=\{x_{i_1,j_1},\ldots,x_{i_n,j_n}\}$, where  
$v_1,\ldots,v_n$ are the corresponding rearrangement of  
$\nu_{i_1,j_1},\ldots,\nu_{i_n,j_n}$, and where $a_1,Q_2,\ldots,Q_n$  
are the corresponding rearrangement of  
$Q_{i_1,j_1},Q_{i_2,j_2},\ldots,Q_{i_n,j_n}$. Thus, if $I$ be the ideal of $R$ generated by $a_{i,j}$ with $(i,j)\in\scrJ_{+,1}$, then the morphism $\Spec \Rscr/\Rscr I\to \Spec R/I$ is a
purely inseparable, finite, flat cover of degree $p^{|\Jscr_{0,0}|}$.

If $D$ is isoclinic, then $\scrJ_{+,2}=\emptyset$, $\kappa(D)=0$, and the condition (C) holds. Thus the morphism $\Spec \Rscr\to \Spec R$ is a finite, flat cover of degree $p^{|\Jscr_{0,0}|}$.
\end{remark}

\begin{remark}{Example}\label{Ex3}
Suppose that $\pi$ is the $r$-cycle $(1\,2\cdots\,r)$. We compute
several invariants introduced above in this Section. We have:
$$\Jscr_{-,1}=\{(i,r)|i\le c\}\cup\{(c,j)|c+1\le j\le
r-1\},$$
$$\Jscr_{+,1}=\{(c+1,j)| j \leq c \}\cup\{(i,1)|c+2\le i\le r\},$$
$$\Jscr_{0,0}=\{(i,j)|r\ge j>i>c\}\cup \{(i,j)|c>i>j\ge 1\}.$$

\begin{itemize}

\item  \noindent Here is a pictorial representation of the sets
for $c=d=4$. The gray diagonal and the arrow remind how the
$r$-cycle $\pi$ acts on the diagram.

\

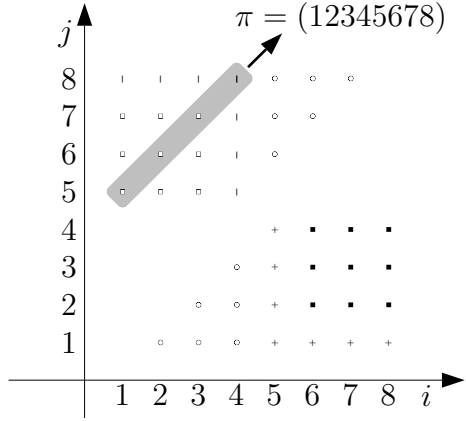
\begin{center}
\begin{pspicture}(-4,-0.5)(9,4.5)
\psset{arrowsize=.2} \psset{arrowlength=1.5} \psset{arrowinset=0}

\rput*[0,0]{45}(2,4){\psframe[linecolor=white,framearc=.5,fillstyle=solid,fillcolor=lightgray](-2.321,.2)(.2,-.2)}
\rput[d](3.37,4.8){$\pi = (12345678)$}

\psline[linewidth=1pt,linearc=.25]{->}(2.15,4.15)(2.6,4.6)

%

\psline[linewidth=0.2pt]{->}(-1,0)(5,0) \rput[d](4.5,-.2){$i$}

\psline[linewidth=0.2pt]{->}(0,-.5)(0,5) \rput[r](-.15,4.6){$j$}

\rput[r](-.1,4){$8$} \psset{dotsize=3pt 0}

\rput[r](-.1,3.5){$7$} \psset{dotsize=3pt 0}

\rput[r](-.1,3){$6$} \psset{dotsize=3pt 0}

\rput[r](-.1,2.5){$5$} \psset{dotsize=3pt 0}

\rput[r](-.1,2){$4$} \psset{dotsize=3pt 0}

\rput[r](-.1,1.5){$3$} \psset{dotsize=3pt 0}

\rput[r](-.1,1){$2$} \psset{dotsize=3pt 0}

\rput[r](-.1,0.5){$1$} \psset{dotsize=3pt 0}


\rput[r](0.6,-0.2){$1$} \psset{dotsize=2pt 0}

\rput[r](1.1,-0.2){$2$} \psset{dotsize=2pt 0}

\rput[r](1.6,-0.2){$3$} \psset{dotsize=2pt 0}

\rput[r](2.1,-0.2){$4$} \psset{dotsize=2pt 0}

\rput[r](2.6,-0.2){$5$} \psset{dotsize=2pt 0}

\rput[r](3.1,-0.2){$6$} \psset{dotsize=2pt 0}

\rput[r](3.6,-0.2){$7$} \psset{dotsize=2pt 0}

\rput[r](4.1,-0.2){$8$} \psset{dotsize=2pt 0}


\multirput(1,0.5)(.5,0){3}{\psdot[dotstyle=o]}
\multirput(1.5,1)(.5,0){2}{\psdot[dotstyle=o]}
\multirput(2,1.5)(.5,0){1}{\psdot[dotstyle=o]}

\multirput(2.5,3)(.5,0){1}{\psdot[dotstyle=o]}
\multirput(2.5,3.5)(.5,0){2}{\psdot[dotstyle=o]}
\multirput(2.5,4)(.5,0){3}{\psdot[dotstyle=o]}

\multirput(2.5,0.5)(.5,0){4}{\psdot[dotstyle=+]}
\multirput(2.5,1)(.5,0){1}{\psdot[dotstyle=+]}
\multirput(2.5,1.5)(.5,0){1}{\psdot[dotstyle=+]}
\multirput(2.5,2)(.5,0){1}{\psdot[dotstyle=+]}

\multirput(2,2.5)(.5,0){1}{\psdot[dotstyle=|]}
\multirput(2,3)(.5,0){1}{\psdot[dotstyle=|]}
\multirput(2,3.5)(.5,0){1}{\psdot[dotstyle=|]}
\multirput(2,4)(.5,0){1}{\psdot[dotstyle=|]}

\multirput(0.5,4)(.5,0){4}{\psdot[dotstyle=|]}

\multirput(3,1)(.5,0){3}{\psdot[dotstyle=square*]}
\multirput(3,1.5)(.5,0){3}{\psdot[dotstyle=square*]}
\multirput(3,2)(.5,0){3}{\psdot[dotstyle=square*]}

\multirput(0.5,2.5)(.5,0){3}{\psdot[dotstyle=square]}
\multirput(0.5,3)(.5,0){3}{\psdot[dotstyle=square]}
\multirput(0.5,3.5)(.5,0){3}{\psdot[dotstyle=square]}

\end{pspicture}
\end{center}

\noindent N.B. We write: $| \in \Jscr_{-,1}$, o $ \in
\Jscr_{0,0}$, $+ \in \Jscr_{+,1}$,
$\square \in \Jscr_{-,2}$, $\blacksquare \in \Jscr_{+,2}$.

\item

\noindent We compute the numbers: $\nu(i,j)$ with
$(i,j)\in\Jscr_{-,1}$. If $(i,j) \in \Jscr_{-,1}$ and $j=r$, then
$\nu(i,j) = c+1-i$; if $(i,j) \in \Jscr_{-,1}$ and $i=c$, then
$\nu(i,j) = r+1-j$. Note that $|\Jscr_{0,0}| = \frac{c(c-1) +
d(d-1)}{2}$ and that $n= |\Jscr_{-,1}| = r-1$. The set of
$\pi$-orders of elements of $\scrJ_{-,1}$ is $\{ 0, 1,2,\ldots, \max\{c,d\} \}$ and therefore $\nu= \max\{c,d\}-1$.

\item Note that the nilpotency bounds of Lemma \ref{NilLem} (b) and
(c) are sharp for $\pi = (12 \cdots r-1 \ r)$. For example, we have  
$\nfr_{0,0}^{c -1} \neq 0$ as by going
down the line $j=i-1$ we have an identity:
$$\bar e_{c,c-1}\bar e_{c-1,c-2}\cdots \bar e_{2,1}=\bar e_{c,1}.$$

\item We compute the invariant $\kappa(\pi)$ explicitly for
$p=2,3$. As $\pi$ is an $r$-cycle, it is easy to see that $\pi$ is
unique up to  conjugation by a permutation of the set $J$ that
leaves invariant the subset $\{1,\ldots,c\}$. Therefore we have
$\kappa(D[p])=\kappa(\pi)$. Suppose that $c\geq 3$ and $d \geq 2$.
The maximal value of $\kappa(\gamma)$ with $\gamma\in
\Gamma_1\cup\Delta_1$ is obtained for certain elements $\gamma \in
\Delta_1$. Like for $\gamma = (c+1,\ c+2,\ \ldots, \ r,\ c,\ 2)$,
we have
$$\kappa(\gamma) = \frac{2}{p} + \frac{1}{p^2} + \frac{1}{p^3} +
\cdots + \frac{1}{p^{d-1}} .$$ For $p=3$, clearly $\kappa(\pi) <
1$. For $p=2$, we have $\kappa(\pi)\geq 1$. Using this one gets that  
if $p=2$ and $c=d\ge 3$, then $D$ is isomorphic to $D\vdual$ and  
the condition $\Cond$ does not hold for $D$.

\item  Suppose $p=2$ and $(c,d)=(3,2)$. Then $\Jscr_{+,1}=\{(4,1),(4,2),(4,3),(5,1)\}$, $\Jscr_{+,2}=\{(5,2),(5,3)\}$, and $\Jscr_{0,0}=\{(2,1),(3,1),(3,2),(4,5)\}$. Thus $\Gamma_1=\{(4,5,2),(4,5,3)\}$, $\Delta_1=\{(4,5,3,2),(5,2,1),(5,3,1)\}$, and $\kappa(\pi)=\kappa(D[p])=1$. The system (14) is  of the form (cf. (15)) 
$$x_{2,1}^4=a_{5,3} x_{4,5}+a_{4,3},\;\;\;\;x_{3,1}^2=a_{5,2} x_{4,5}+a_{5,3} x_{4,5}x_{2,1}^2+a_{4,2}+a_{4,3}x_{2,1}^2$$
$$x_{4,5}^2=a_{5,2} x_{2,1}+a_{5,3} x_{3,1}+a_{5,1}, \;\;\;\;x_{4,1}=a_{4,1}+a_{4,2}x_{2,1}+a_{4,3}x_{3,1}.$$
Eliminating $x_{4,1}$ and substituting $x_{3,1}=t^{7/4}$, Lemma \ref{GrobLem} implies that this system defines a finite morphism  
(i.e., $\zeta$ is a finite
morphism despite the fact that $\kappa(\pi)=\kappa(D[p])=1$). This
points out that the general weights $\mu_{\ell}$ used in the proof of
Theorem \ref{FinThm} can occasionally be replaced by other weights
which can still make the proof of Theorem \ref{FinThm} work even if $\kappa(D[p])\ge 1$.
\end{itemize}

\end{remark}

\end{segment}


\section{Proof of Theorem~\ref{MainThm}}\label{MainProof}

In this Section, we prove Theorem~\ref{MainThm} (see Subsection
\eqref{ProofofMain}).

\begin{segment}{Exact subgroups}{ExSubGrp}
We first recall from \cite{CPS} several properties of exact subgroups  
of an affine,
smooth, connected group $\Gscr$ over a field which are needed in Subsection
\eqref{ProofofMain}. A smooth subgroup scheme $\Fscr$ of $\Gscr$ is called
\emph{exact} if the induction of rational $\Fscr$-modules to rational
$\Gscr$-modules preserves short exact sequences.

\begin{proclamation}{Theorem}\label{Thm51}
The following six properties hold:
\begin{definitionlist}
\item The smooth subgroup scheme $\Fscr$ of $\Gscr$ is exact if
and only if the quotient variety $\Gscr/\Fscr$ is affine. \item
The smooth subgroup scheme $\Fscr$ of $\Gscr$ is exact if and only
if its identity component $\Fscr^0$ is exact. \item If $\Gscr$ is
unipotent, then each smooth subgroup of it is exact. \item If
$\Kscr$ is a smooth subgroup of $\Fscr$ which is exact in $\Fscr$
and if $\Fscr$ is exact in $\Gscr$, then $\Kscr$ is exact in
$\Gscr$. \item Let $\Nscr$ be a normal, smooth, connected subgroup
of $\Gscr$. Let $\tilde\Fscr:=\Im(\Fscr\to\Gscr/\Nscr)$; it is a
smooth subgroup of $\Gscr/\Nscr$. Then $\tilde\Fscr$ is exact in
$\Gscr/\Nscr$ if and only if the reduced group
$(\Nscr\Fscr)_{\text{red}}$ is exact in $\Gscr$. \item If $\Fscr$
is exact in $\Gscr$, then it is exact in every other smooth,
connected subgroup $\Kscr$ of $\Gscr$ that contains $\Fscr$.
\end{definitionlist}
\end{proclamation}

\begin{proof} Part (a) is the main result of \cite{CPS} (cf.  
\cite{CPS}, Theorem 4.3).
The natural morphism $\Gscr/\Fscr^0\to\Gscr/\Fscr$ is finite and
surjective. Thus $\Gscr/\Fscr^0$ is an affine variety if and only
if $\Gscr/\Fscr$ is an affine variety, cf. Chevalley's theorem of
\cite{EGAII}, Th\'eor\`eme (6.7.1). From this and (a), we get that
(b) holds. Parts (c) and (d) are implied by \cite{CPS}, Corollary
2.2 and \cite{CPS}, Proposition 2.3 (respectively).  Part (e)
follows from (a) once we remark that $(\Gscr/\Nscr)/\tilde\Fscr$
is isomorphic to $\Gscr/(\Nscr\Fscr)_{\text{red}}$. Part (f)
follows from \cite{CPS}, Proposition 2.4.
\end{proof}

\begin{proclamation}{Lemma}\label{Lem52}
We assume that the orbit $\Oscr_1$ is an affine variety over $k$. Then  
for each $m\in\Bbb N^*$, the orbit $\Oscr_m$ is an affine variety over $k$.
\end{proclamation}

\begin{proof}
As we have a finite, surjective morphism
$\Hscr_1/\Cscr_1^0\to\Hscr_1/\Sscr_1\arrowsim \Oscr_1$, from
our hypothesis we get that the quotient variety $\Hscr_1/\Cscr_1^0$ is  
affine. Thus $\Cscr_1^0$ is exact in $\Hscr_1$, cf.  Theorem  
\ref{Thm51} (a).

Let $\Cscr_{1,m}^0$ be the pullback of $\Cscr_1^0$ via
the reduction epimorphism $\Hscr_m\twoheadrightarrow\Hscr_1$. As
$\Ker(\Hscr_m\twoheadrightarrow\Hscr_1)$ is a unipotent group (cf.
\cite{Va2}, Lemma 2.2.3) which is smooth and connected and as
$\Cscr_1^0$ is a unipotent, smooth, connected group, we get that
the group $\Cscr_{1,m}^0$ is a unipotent, smooth, connected group.
As $\Cscr_1^0$ is exact in $\Hscr_1$, from Theorem \ref{Thm51} (e)
we get that $\Cscr_{1,m}^0$ is exact in $\Hscr_m$. The subgroup
$\Cscr_m^0$ of $\Cscr_{1,m}^0$ is exact in $\Cscr_{1,m}^0$, cf.
Theorem \ref{Thm51} (c). From the last two sentences and Theorem
\ref{Thm51} (d), we get that $\Cscr_m^0$ is exact in $\Hscr_m$.
Thus the quotient variety $\Hscr_m/\Cscr_m^0$ is affine, cf.
Theorem \ref{Thm51} (a).

As we have a finite, surjective morphism
$\Hscr_m/\Cscr_m^0\to\Hscr_m/\Sscr_m\arrowsim \Oscr_m$, we get
from Chevalley's theorem that $\Oscr_m$ is an affine variety over
$k$.
\end{proof}
\end{segment}

\begin{segment}{Proof of Theorem~\ref{MainThm}}{ProofofMain}
Based on Subsection\eqref{devissage}, we can assume that we are in the Essential Situation \ref{EssSit}. In
particular, $D_m=D[p^m]$ and $k$ is algebraically closed. We
recall that we are assuming that the condition $\Cond$ holds for
$D_m$ (see the Introduction) and thus also for $D$. As
$\grs_{D}(m)=\Ascr_D(\Escr)$ is equal to $\Ascr_{D\vdual}(\Escr\vdual)$, by
replacing if needed $(D,\Escr)$ with
$(D\vdual,\Escr\vdual)$ we can assume that
$\kappa(D)<1$. Therefore we can choose the permutation $\pi$ of
Subsection \eqref{Nilp} such that $\kappa(\pi)<1$. We also recall
that $\tilde\Cscr_{1}^0$ is a connected, smooth, unipotent
subgroup of $\Iscr_{1}$. As $\kappa(\pi)<1$, from
Theorem~\ref{FinThm} we get that we have a natural finite, flat
morphism
$\Vscr_2\times_k\Wscr_{0,k}\twoheadrightarrow\tilde\Cscr_{1}^0\backslash\Iscr_{1}$
of degree $p^{|\Jscr_{0,0}|}$ (in the smooth, faithfully flat
topology, this morphism looks like the morphism $\zeta$). From
this and Chevalley's theorem, we get that the quotient variety
$\tilde\Cscr_{1}^0\backslash\Iscr_{1}$ is affine. As the quotient
varieties $\tilde\Cscr_{1}^0\backslash\Iscr_{1}$ and
$\Iscr_1/\tilde\Cscr_1^0$ are isomorphic, we get that
$\Iscr_1/\tilde\Cscr_1^0$ is affine. Thus $\tilde\Cscr_{1}^0$ is
an exact subgroup of $\Iscr_{1}$, cf. Theorem \ref{Thm51} (a).
Therefore the unipotent group
$(\Wscr_{-,k}\Cscr_{1}^0)_{\text{red}}$ is an exact subgroup
scheme of $\Hscr_{1}$, cf. Theorem \ref{Thm51} (e). From this and
Theorem \ref{Thm51} (c) and (d), we get that $\Cscr_{1}^0$ is an
exact subgroup of $\Hscr_{1}$. As we have a finite, surjective
morphism $\Hscr_1/\Cscr_1^0\to\Hscr_1/\Sscr_1\arrowsim \Oscr_1$,
from Theorem \ref{Thm51} and Chevalley's theorem, we get that
$\Oscr_1$ is affine. From this and Lemma \ref{Lem52}, we get that:

\begin{proclamation}{Corollary}\label{AffOrb}
We assume that the condition $\Cond$ holds for $D$. Then for each $m\in\Bbb N^*$, the orbit  
$\Oscr_m$ is an
affine variety over $k$.
\end{proclamation}

We can complete the proof of Theorem~\ref{MainThm} as follows.
Locally in the Zariski topology of $\Ascr$, we can write $\Ascr=\Spec A$ and we can assume that the evaluation of the Dieudonn\'e crystal $\Bbb D(\Escr)$ of $\Escr$ at the thickening $\Ascr\hookrightarrow \Spec W(A)$ is, when viewed without connection and Verschiebung map, a pair of the form $(M\otimes_{W(k)} W(A),g_A(\phi\otimes\sigma_A))$ with $g_A\in\pmb{GL}_M(W(A))$. The reduction $g_A[m]\in \pmb{GL}_M(W_m(A))=\Dscr_m(A)$ of $g_A$ is a morphism $\varpi:\Ascr\to \Dscr_m$. The locally closed
embedding $\grs_D(m)\hookrightarrow \Ascr$ is the pullback of the
locally closed embedding $\Oscr_m\hookrightarrow\Dscr_m$ via $\varpi$ (cf.
\cite{Va2}, Subsubsection 3.1.1). 

From the previous paragraph and Corollary
\ref{AffOrb}, we get that in general $\grs_D(m)$ is an affine $\Ascr$-scheme.
This ends the proof of Theorem~\ref{MainThm}.\endproof

\begin{remark}{Remark}\label{Aut}
Under the natural identification $\Iscr_1=\Wscr_{+0,k}$, the group
$\tilde\Cscr_1^0$ gets identified with the crystalline realization $\pmb{Aut}(D[p])_{\text{crys},\text{red}}^0$ of the identity component $\pmb{Aut}(D[p])_{\text{red}}^0$ of the reduced group scheme of automorphisms of $D[p]$ (cf. \cite{Va2}, Theorem 2.4 (b)). Thus the quotient variety $\Wscr_{+0,k}/\pmb{Aut}(D[p])_{\text{crys},\text{red}}^0$
is affine. It is easy to see that $\nfr_{0,0}\oplus\nfr_{+,1}$ is the  
$k$-span of elements $g-1_M[1]$ with $g\in  
\pmb{Aut}(D[p])_{\text{crys},\text{red}}^0(k)=\pmb{Aut}(D[p])_{\text{red}}^0(k)$. Therefore the number  
$|\Jscr_{0,0}|$ is an invariant of the isomorphism class of $D[p]$.
\end{remark}
\end{segment}


\section{Level $m$ stratifications for quasi Shimura $p$-varieties of  
Hodge type}\label{StratHodge}

In this Section, we consider the relative version of the action
$\Bbb T_m$. We prove Proposition \ref{GenAffOrb}, an analogue of
Corollary \ref{AffOrb}, which will imply Theorem \ref{SiegelThm}
and its analogues for general level $m$ stratifications.

\begin{segment}{Orbits of relative group actions}{PureLevelm}
Let $G$ be a smooth, closed subgroup scheme of $\pmb{GL}_M$ such that  
its generic
fibre $G_{B(k)}$ is connected. Thus the scheme $G$ is integral. Until  
the end we
will assume that the following two axioms hold for the triple $(M,\phi,G)$:

\medskip
\begin{romanlist}
\item
The Lie subalgebra $\Lie(G_{B(k)})$ of $G_{B(k)}$ is stable under  
$\phi$ i.e., we
have $\break\phi(\Lie(G_{B(k)}))=\Lie(G_{B(k)})$.
\item
There exist a direct sum decomposition $M=F^1\oplus F^0$ and a smooth, closed
subgroup scheme $G_1$ of $\pmb{GL}_M$ such that the following four 
properties hold:
\begin{definitionlist}
\item
The kernel of the reduction modulo $p$ of $\phi$ is $F^1/pF^1$.
\item
The cocharacter $\mu:\Bbb G_m\to\pmb{GL}_M$ which fixes $F^0$ and acts via the
inverse of the identical character of $\Bbb G_m$ on $F^1$, factors  
through $G_1$.
\item
The group scheme $G_1$ contains $G$ in such a way that we have a
short exact sequence $1\to G\to G_1\to\Bbb G_m^u\to 1$ with $u\in\{0,1\}$.
\item
If $u=1$ (i.e., if $G_1\neq G$), then the homomorphism $\mu:\Bbb G_m\to  
G_1$ defined
by $\mu$ (cf. properties (b) and (c)) is a splitting of the short exact sequence of the
property (c).
\end{definitionlist}
\end{romanlist}

\interbreak

\noindent We will use the notations of Section \ref{GroupActions}
for the direct sum decomposition $M=F^1\oplus F^0$ of the axiom
(ii). Due to the properties (ii.b) and (ii.c) we have a direct sum
decomposition into $W(k)$-modules
$$\Lie(G)=\oplus_{i=-1}^1\tilde F^i(\Lie(G))\leqno (17)$$
such that $\mu$ acts via inner conjugation on $\tilde
F^i(\Lie(G))$ as the $-i$-th power of the identical character of
$\Bbb G_m$. Let $e_+$, $e_0$, and $e_-$ be the ranks of $\tilde
F^1(\Lie(G))$, $\tilde F^0(\Lie(G))$, and $\tilde F^{-1}(\Lie(G))$
(respectively). Let $d_G:=\dim(G_k)=\dim(G_{B(k)})$. Due to (17),
we have
$$d_G=e_++e_0+e_-.$$
We consider the following four closed subgroup schemes
$\Wscr_+^G:=\Wscr_+\cap G$, $\Wscr_0^G:=\Wscr_0\cap G$,  
$\Wscr_-^G:=\Wscr_-\cap G$,
and $\Wscr_{+0}^G:=\Wscr_{+0}\cap G$ of $G$. Let
$$\Hscr^G:=\Wscr_+^G\times_{W(k)} \Wscr_0^G\times_{W(k)} \Wscr_-^G;$$
it is a closed subscheme of $\Hscr$ such that $\Hscr^G_{W_m(k)}$ is a closed
subgroup subscheme of $\Hscr_{W_m(k)}=\tilde\Hscr_{W_m(k)}$ (we recall  
that we view
the isomorphism  
$\Pscr_{W_m(k)}:\Hscr_{W_m(k)}\arrowsim\tilde\Hscr_{W_m(k)}$ of
$\Spec W_m(k)$-schemes as a natural identification).

\noindent The group schemes $\Wscr_+^G$ and $\Wscr_-^G$ over
$\Spec W(k)$ are isomorphic to $\Bbb G_a^{e_+}$ and $\Bbb G_a^{e_-}$
(respectively). More precisely, if $R$ is a commutative
$W(k)$-algebra, then we have identities $\Wscr_+^G(R)=1_{M\otimes_{W(k)}
R}+\tilde F^1(\Lie(G))\otimes_{W(k)}
R$ and $\Wscr_-^G(R)=1_{M\otimes_{W(k)} R}+\tilde
F^{-1}(\Lie(G))\otimes_{W(k)} R.$ \indent Let
$\Wscr_0^{G_1}:=\Wscr_0\cap G_1$. The group schemes
$\Wscr_0^{G_1}$ and $\Wscr_0^G$ are smooth (see \cite{Va2},
Subsubsection 4.1.1).

\noindent The Lie algebras of $W_+^G$, $W_0^G$, and $W_-^G$ are
$\tilde F^1(\Lie(G))$, $\tilde F^0(\Lie(G))$, and $\tilde
F^{-1}(\Lie(G))$ (respectively). This implies that the relative
dimension of $\Wscr_0^G$ is $e_0$. The smooth, affine scheme
$\Hscr^G$ has relative dimension $d_G$ over $\Spec W(k)$. The
natural product morphism $\Pscr_0^G:\Hscr^G\to G$ is induced
naturally by the open embedding
$\Pscr_0:\Hscr\hookrightarrow\pmb{GL}_M$ and therefore it is also
an open embedding.  Let $\Pscr_-^G:=1_{\Wscr_+^G}\times
1_{\Wscr_0^G}\times p1_{\Wscr_-^G}:\Hscr^G\to\Hscr^G$. The
composite morphism
$\Pscr_{0-}^G:=\Pscr_0^G\circ\Pscr_-^G:\Hscr^G\to G$ has the
property that its reduction
$\Pscr_{0-,W_m(k)}^G:\Hscr^G_{W_m(k)}\to G_{W_m(k)}$ modulo $p^m$
is a homomorphism of affine group schemes over $\Spec W_m(k)$
which is a restriction of the homomorphism
$\Pscr_{0-,W_m(k)}:\Hscr_{W_m(k)}\to\pmb{GL}_{M/p^mM}$ (see
Section \ref{GroupActions} for $\Pscr_{0-}$).

\noindent Let $\tilde\Hscr^G$ be the dilatation of $G$ centered on
the smooth subgroup $\Wscr^G_{+0,k}$ of $G_k$; it is a smooth,
closed subgroup scheme of $\tilde\Hscr$. As in Section
\ref{GroupActions}, we argue that we have a natural morphism
$\Pscr^G:\Hscr^G\to\tilde\Hscr^G$ of $\Spec W(k)$-schemes which
gives birth to an isomorphism
$\Pscr_{W_m(k)}^G:\Hscr^G_{W_m(k)}\arrowsim\tilde\Hscr^G_{W_m(k)}$
of $\Spec W_m(k)$-schemes, to be viewed as an identification. The
group schemes structures on $\Hscr_{W_m(k)}^G$ induced via the
identification $\Pscr^G_{W_m(k)}$ or via the identification of
$\Hscr_{W_m(k)}^G$ with a closed subgroup scheme of
$\Hscr_{W_m(k)}$, are equal.

\noindent Let $\Hscr_m^G:=\Bbb W_m(\Hscr^G)=\Bbb W_m(\tilde\Hscr^G)$;
it is a smooth, affine group over $k$ of dimension $md_G$ which is
connected if and only if $\Hscr_k^G$ (equivalently $\Wscr_{0,k}^G$) is  
connected (cf. Section \ref{GroupActions}).
Let $\Dscr_m^G:=\Bbb W_m(G)$; it is a smooth, affine $k$-scheme of
dimension $md_G$ which is connected if and only if $G_k$ is
connected (cf. Section \ref{GroupActions}). There exists a unique action
$$\Bbb T_m^G:\Hscr_m^G\times_k\Dscr_m^G\to\Dscr_m^G$$
which is the natural restriction of the action $\Bbb T_m$ of Section
\ref{GroupActions}.

\noindent
Let $\Oscr_m^G$ be the orbit of $1_M[m]\in\Dscr_m^G(k)$ under
the action $\Bbb T_m^G$. Let $\Sscr_m^G$ be the stabilizer subgroup  
scheme of the
point $1_M[m]\in\Dscr_m^G(k)$ under the action $\Bbb T_m^G$; we have
$\Sscr_m^G:=\Sscr_m\cap\Hscr_m^G$. Let $\Cscr_m^G$ be the reduced group of
$\Sscr_m^G$ and let $\Cscr_m^{0G}$ be the identity component of $\Cscr_m^G$.

\begin{proclamation}{Proposition}\label{GenAffOrb}
If the condition $\Cond$ holds for $D$, then the orbit $\Oscr_m^G$
of $1_M[m]\in\Dscr_m^G(k)$ under the action $\Bbb T_m^G$ is a
smooth, affine $k$-scheme. Therefore, if either $p=3$ and
$\min\{c,d\} \leq 6$ or $p\geq 5$, then each orbit of the action $\Bbb T_m^G$ is
a smooth, affine $k$-scheme.
\end{proclamation}

\begin{proof}
As $\Cscr_m^0$ is a unipotent group, $\Cscr_m^{0G}$ is exact in
$\Cscr_m^0$ (cf. Theorem \ref{Thm51} (c)). From this and the fact
that $\Cscr_m^0$ is exact in $\Hscr_m$ (cf. Subsection
\eqref{ProofofMain}), we get that $\Cscr_m^{0G}$ is exact in
$\Hscr_m$ (cf. Theorem \ref{Thm51} (d)) and thus also in
$\Hscr_m^G$ (cf. Theorem \ref{Thm51} (f)). As in Subsection
\ref{ProofofMain}, we argue that this implies that  $\Oscr_m^G$ is
a smooth, affine $k$-scheme.
\end{proof}
\begin{remark}{Example}\label{Ex4}
Let $a\geq 2$ be an integer. We assume that we have a direct
sum decomposition $M=\oplus_{s=1}^a M_s$ such that $F^1\subseteq M_1$ and for all
$s\in\{1,\ldots,a\}$ we have $\phi(M_s)\subseteq M_{s+1}$ (here
$M_{a+1}:=M_1$). We identify naturally $W(\Bbb F_{p^a})$ with a
$\Bbb Z_p$-subalgebra of $\End(M,\phi)$ (equivalently of $\End(D)$)
that acts on each $M_s$ via scalar endomorphisms. From Example
\ref{Ex1}, we get that the condition $\Cond$ holds for $D$. Thus $\Oscr_m^G$ is a smooth, affine $k$-scheme, cf. Proposition \ref{GenAffOrb}. Thus, if $G$ is a closed subgroup scheme of $\prod_{s=1}^a \pmb{GL}_{M_s}$, then each orbit of $\Bbb T_m^G$ is a smooth, affine $k$-scheme. We emphasize
that these hold for all primes $p$.
\end{remark}
\end{segment}

\begin{segment}{Quasi Shimura $p$-varieties of Hodge type}{QuasiShimura}
In this Subsection, we assume that $c=d$ and we use only $d$. We also  
assume that the condition $\Cond$ holds for
each $p$-divisible group over $k$ which admits principal
quasi-polarizations and has dimension $d$ (e.g., if $p \geq 5$).  
Suppose that $D$
has a principal quasi-polarization $\lambda$. Let $\psi:M\times
M\to W(k)$ be the perfect, alternating form on $M$ induced
naturally by $\lambda$; for $x,y\in M$, we have
$\psi(\phi(x),\phi(y))=p\sigma(\psi(x,y))$. Suppose that $G$ is a
closed subgroup scheme of $\pmb{Sp}(M,\psi)$. We recall that axioms  
(i) and (ii) of Subsection \eqref{PureLevelm} hold for the
triple $(M,\phi,G)$. As $\mu:\Bbb G_m\to G_1$ cannot factor through
$\pmb{Sp}(M,\psi)$, we have $u=1$ (i.e., we have a short exact
sequence $1\to G\to G_1\to\Bbb G_m\to 1$). Let
$$\Fscr:=\{(M,g\phi,\psi,G)|g\in G(W(k))\};$$
it is a {\it family} of principally quasi-polarized Dieudonn\'e
modules with a group over $k$. Let $\Ascr_{d,1,N}$ be as in  
Introduction. 

\noindent
Let $\Mscr$ be a {\it quasi Shimura
$p$-variety of Hodge type relative to $\Fscr$} in the sense of
\cite{Va2}, Definition 4.2.1. Thus $\Mscr$ is a smooth $k$-scheme which is
equidimensional and which is equipped with a morphism
$\Mscr \to \Ascr_{d,1,N,k}$ that induces $k$-epimorphisms at the
level of complete, local rings of residue field $k$ (i.e., it is a
formal closed embedding at all $k$-valued points) and that
satisfies an extra axiom that involves $\Fscr$ (see \cite{Va2}, Axiom
4.2.1 (iii)). This extra axiom implies that $\Mscr$ has a level
$m$ stratification $\Sscr^G(m)$ in the sense of \cite{Va1}, Definition 2.1.1.
For each algebraically closed field $K$ that contains $k$, we have
a set $\Sscr_m^G(K)$ of reduced locally closed subschemes of
$\Mscr_K$ which are regular, equidimensional and which locally in
the \'etale topology of $\Mscr_K$ are pullbacks of locally closed  
embeddings of the form $\Oscr^G\hookrightarrow
\Hscr_{m,K}^G$ for some orbit $\Oscr^G$ of the extension of the action
$\Bbb T_m^G$ to $K$ (cf. \cite{Va2}, Subsubsection 4.2.3). From this and
the analogue of the Proposition \ref{GenAffOrb} over $K$, we get that each
element of $\Sscr_m^G(K)$ is pure in $\Mscr_K$. Thus we
get that:

\begin{proclamation}{Theorem}\label{PurityShimura}
The level $m$ stratification $\Sscr_m^G$
of $\Mscr$ is pure. In other words, for each
algebraically closed field $K$ that contains $k$, all the elements
of $\Sscr_m^G(K)$ (i.e., all strata of $\Sscr_m^G$ which are
subschemes of $\Mscr_K$) are affine $\Mscr_K$-schemes
(equivalently, are affine $\Mscr$-schemes).
\end{proclamation}

\begin{proclamation}{Corollary}\label{CodimShimura}
Let $K$ be an algebraically closed field extension of $k$. Let
$\grs\in\Sscr_m^G(K)$. Let $\bar\grs$ be the Zariski
closure of $\grs$ in $\Mscr_K$. Then the reduced scheme underlying
$\bar\grs\setminus\grs$ is either empty or of pure codimension $1$
in $\bar\grs$.
\end{proclamation}

\begin{remark}{Example}\label{Ex5}
Suppose that $G=\pmb{\text{Sp}}(M,\psi_M)$ and  
$G_1=\pmb{\text{GSp}}(M,\psi_M)$.
Then $\Ascr_{d,1,N,k}$ is a quasi Shimura $p$-variety of Hodge type  
relative to
$\Fscr$, cf. \cite{Va2}, Example 4.5. The strata of $\Sscr_m^G(k)$ are  
of the form
$\grs_{D,\lambda}(m)$, cf. loc. cit. Thus Theorem \ref{SiegelThm} is a  
particular
case of Theorem \ref{PurityShimura}.
\end{remark}

\begin{remark}{Remark}\label{Primes23}
One can use Example \ref{Ex4} to get plenty of examples of level $m$  
stratifications
in characteristic $2$ or $3$ which are pure.
\end{remark}
\end{segment}

\begin{section}{On principal purity for stratifications}\label{StrongPurity}

In Definition~\ref{DefPrincPure}, we introduced the stronger notion
for a subscheme to be principally pure (with respect to some
Grothendieck topology). In this section, we investigate the
question whether the $p$-rank stratification is principally pure.
This is connected to the existence of generalized Hasse--Witt
invariants i.e., to the local existence of functions $f$ such that
a stratum is defined in its scheme-theoretic closure as the locus where $f$ is
invertible.

In \cite{Ito}, Theorem on p. 1567, T. It\=o proves the existence of
generalized Hasse--Witt invariants for $p$-rank strata of the
special fibre (at a split prime $p$) of a good integral model
of a unitary Shimura variety $\text{Sh}(\Gscr,\Xscr)$ which is of
PEL type and which over $\Bbb R$ involves a derived group
$\Gscr_{\Bbb R}^{\text{der}}=\pmb{SU}(n-1,1) \times \pmb{SU}(n,0) \times \cdots \times\pmb{SU}(n,0)$. In this case, the $p$-rank strata coincide with
the level $m$ strata for all $m$. He derives the corollary that
the $p$-rank strata are affine. From the affineness
result and the weak Lefschetz theorem, he obtains that the number
of connected components of all Zariski closures of positive
dimensional $p$-rank strata are equal.

In general, however, principal purity (for the \'etale topology)
does not hold for $p$-rank strata, even over regular schemes, as
the following example shows.

\begin{remark}{Example}\label{Ex6}
We set $c = d = 2$. Let $\Rscr:=k\dlbrack
x_1,x_2,x_3,x_4\drbrack$ and $R:=\Rscr/(x_1x_4 - x_2x_3)$. Let $S :=
\Spec \Rscr$ and $T:=\Spec R$. Let $g$ be the $\Rscr$-linear  
automorphism of $N:=\Rscr^4$ associated
to an invertible block matrix
\begin{displaymath}
\mathbf A =
\left(\begin{array}{cc}
A_1 & A_2\\
A_3 & A_4
\end{array}
\right)
\end{displaymath}
formed by $2\times 2$ blocks, where
\begin{displaymath}
A_1 =
\left(\begin{array}{cc}
x_1 & x_2\\
x_3 & x_4
\end{array}
\right).
\end{displaymath}
Let $\{ \bar e_1,\bar e_2,\bar e_3,\bar
e_4\}$ be the standard basis for $N$. Let $\phi_N:N\otimes_{\Rscr, \sigma_{\Rscr}} \Rscr\to N$ be the
$\Rscr$-linear map that takes the quadruple
$(\bar e_1\otimes 1,\bar e_2\otimes 1, \bar e_3\otimes 1,\bar  
e_4\otimes 1)$ to the
quadruple
$(g(\bar e_1),g(\bar e_2),0,0)$. Let $\vartheta_N:N\to  
N\otimes_{\Rscr, \sigma_{\Rscr}}\Rscr$
be the $\Rscr$-linear map that takes the quadruple
$(g(\bar e_1),g(\bar e_2),g(\bar e_3),g(\bar e_4))$ to the quadruple  
$(0,0,\bar
e_3\otimes 1,\bar e_4\otimes 1)$. It is easy to see that the triple
$(N,\phi_N,\vartheta_N)$ is the Dieudonn\'e crystalline functor  
(viewed without
connection) of a $BT_1$ over $S$, to be denoted as $B$.

\noindent 
One computes the $p$-rank of $B$ at a geometric point $x$ of $S$ as follows. Let $N_1$ be the direct summand of $N$ generated by $\bar e_3$ and $\bar e_4$. The kernel of $\phi_N$ is $N_1\otimes_{\Rscr, \sigma_{\Rscr}}\Rscr$. Let $\bar\phi_N: (N/N_1)\otimes_{\Rscr, \sigma_{\Rscr}}\Rscr\to N/N_1$ be the $\Rscr$-linear map induced naturally by $\phi_N$. We view $\{\bar e_1,\bar e_2\}$ as an $\Rscr$-basis for $N/N_1$. One looks for solutions of the equation
$\bar\phi_N((z_1\bar e_1+z_2 \bar e_2)\otimes 1)=z_1 \bar e_1+z_2 \bar e_2$. One comes across the following system $\Yscr$ of two equations
$$ z_1 = x_1 z_1^p + x_2 z_2^p,\;\;\;\;z_2 = x_3 z_1^p + x_4 z_2^p.$$
The $p$-rank of $B$ at $x$ is the dimension of the $\Bbb F_p$-vector space of solutions of $\Yscr$ at $x$. Therefore the $p$-rank stratification of $S$ has an open dense
stratum $S\setminus T$ (of $p$-rank $2$) and has one stratum $T\setminus Y$ of
codimension $1$ in $S$ (of $p$-rank $1$). Here $Y$ is the smallest (thus reduced) closed subscheme
of $T$ with the property that $\Yscr$ defines a scheme over $T$ which is an \'etale cover of degree $p$ above $T\setminus Y$. 

We show that the assumption that $Y$ is the zero locus of a single
function $f\in R$ leads to a contradiction. To solve $\Yscr$ over $T$, we remark that since we have $x_1x_4=x_2x_3$ in $R$, we also have $x_3z_1=x_1z_2$. Thus to solve $\Yscr$ over $T':=\Spec R[{1\over x_1}]$, we can substitute $z_2=x_1^{-1}x_3z_1$ into the first equation of $\Yscr$ and get the equation $z_1=(x_1+x_1^{-p}x_2x_3^p)z_1^p$. One concludes that $Y\cap T'$ is the zero locus of the function $x_1+x_1^{-p}x_2x_3^p=x_1^{1-p}(x_1^p+x_4x_3^{p-1})$ in $T'$. One checks that $T'$ is regular and that $\Spec R[{1\over {x_1}}]/(x_1^p+x_4x_3^{p-1})$ is regular and irreducible (irreducibility can be checked starting from the fact that $k[x_1,x_2,x_3,x_4]/(x_1x_4-x_2x_3,x_1^p+x_4x_3^{p-1})[{1\over {x_1}}]$ is isomorphic to $k[x_1,x_3][{1\over {x_1x_3}}]$). Moreover, each unit in $R[{1\over {x_1}}]$ is a unit of $R$ times an integral power of $x_1$. From the last two sentences we get that the image of $f$ in $R[{1\over {x_1}}]$ is equal to $(x_1^p+x_4x_3^{p-1})^vx_1^tu$, where $v$ is a positive integer, where $t$ is an integer, and where $u$ is a unit in $R$. It is easy to see that $t\ge 0$ and therefore the equations $x_1=x_3=x_2=0$ define a closed subscheme of $Y$. But over the locus defined by $x_4\neq 0$ and $x_1=x_2=x_3=0$, the $p$-rank of $B$ is $1$, and this is a contradiction to the $p$-rank being $0$ there.
Therefore, principal purity fails for the $p$-rank $1$ stratum
$T\setminus Y$ of $S$. 
\end{remark}

\begin{segment}{Proof of Proposition~\ref{NoPrincPure}}{Proofof18}
Let $S$ and $B$ be as in Example
\ref{Ex6}. From Example \ref{Ex6}, we get that the Proposition \ref{NoPrincPure} holds for $c=d=2$ and $s=1$. Using
direct sums of $B$ and of constant $BT_1$'s over $S$, one easily gets that
the Proposition \ref{NoPrincPure} holds for all $c,d\geq 2$ and $s\in\{1,\ldots,c-1\}$.\endproof
\end{segment}

\end{section}


\bigskip\noindent
{\bf Acknowledgments.} The first author was supported by a
postdoctoral research fellowship of the FQRNT while enjoying the
hospitality of the Institut de math\'ematiques de Jussieu
(Universit\'e Paris $7$); he thanks O. Brinon for kindly sharing
his TeX expertise. The second author would like to thank
Binghamton University for good working conditions; he was partially supported by the NSF grant \#0900967.
\bigskip


\bigskip\bigskip

\hbox{Marc--Hubert Nicole}
\hbox{Email: nicole@iml.univ-mrs.fr}
\hbox{Address: Institut math\'ematique de Luminy,}
\hbox{Campus de Luminy, Case 907,}
\hbox{13288 Marseille cedex 9, France.}

\bigskip
\hbox{Adrian Vasiu}
\hbox{Email: adrian@math.binghamton.edu}

  \hbox{Address: Department of Mathematical Sciences,
Binghamton University,} \hbox{Binghamton, NY 13902-6000, U.S.A.}

\bigskip
\hbox{Torsten Wedhorn}
\hbox{Email: wedhorn@math.uni-paderborn.de} \hbox{Address:
Universit\"at Paderborn, Institut f\"ur Mathematik,}
\hbox{Warburger Str. 100, 33098 Paderborn, Germany.}

\bigskip
\noindent
Final version to appear in Ann. Sci. \'Ec. Norm. Sup.


\begin{thebibliography}{EGAIII}
\bibitem[BLR]{BLR} S.~Bosch, W.~L\"utkebohmert, M.~Raynaud, {\em  
N\'eron Models},
Ergebnisse der Mathematik und ihrer Grenzgebiete. 3. Folge, {\bf 21},
Springer-Verlag, Berlin, 1990.\\[-20pt]
\bibitem[Bo]{Bo} P. Boyer, Monodromie du faisceau pervers des cycles \'evanescents de quelques vari\'et\'es de Shimura
simples, \emph{Invent. Math.} {\bf 177}  (2009),  no. 2, 239--280.\\[-20pt]
\bibitem[CPS]{CPS} E.~Cline, B.~Parshall, L.~Scott, Induced modules and affine
quotients, {\em Math. Annalen} {\bf 230} (1977), no. 1, 1--14.\\[-20pt]
\bibitem[De]{De} P.~Deligne, Vari\'et\'es de Shimura: interpr\'etation  
modulaire, et
techniques de construction de mod\`eles canoniques, in {\em  
Automorphic forms, representations, and $L$-functions} (Oregon State  
Univ., Corvallis, OR, 1977), Part 2,  247--289, Proc. Sympos. Pure  
Math., Vol. {\bf 33}, Amer. Math. Soc., Providence, RI, 1979.\\[-20pt]
\bibitem[dJO]{dJO} J.~de Jong, F.~Oort, Purity of the stratification  
by Newton polygons, {\em J. Amer. Math. Soc.} {\bf 13} (2000), no. 1,  
209--241.\\[-20pt]
\bibitem[EGAI]{EGAI} A.~Grothendieck, J.~Dieudonn\'e, {\em \'El\'ements de
G\'eom\'etrie Alg\'ebrique, I. Le langage des sch\'emas}, Grundlehren  
der Mathematik {\bf 166}, Springer-Verlag, 1971.\\[-20pt]
\bibitem[EGAII]{EGAII} A.~Grothendieck, J.~Dieudonn\'e, {\em \'El\'ements de
G\'eom\'etrie Alg\'ebrique, II. \'Etude globale \'el\'ementaire de
quelques classes de morphismes}, Inst. Hautes \'Etudes Sci. Publ.  
Math. No. {\bf 8}, 1961.\\[-20pt]
\bibitem[EGAIV]{EGAIV} A.~Grothendieck, J.~Dieudonn\'e, {\em \'El\'ements de
G\'eom\'etrie Alg\'ebrique, IV. \'Etude locale des sch\'emas et des  
morphismes de
sch\'emas, Seconde partie}, Inst. Hautes \'Etudes Sci. Publ. Math. No.
{\bf 24}, 1965.\\[-20pt]
\bibitem[Go]{Go} E. Z. Goren, Hasse invariants for Hilbert
modular varieties, {\em Israel J. Math.} {\bf 122} (2001),
157--174.\\[-20pt]
\bibitem[Ha]{Ha} G. Harder, \"Uber die Galoiskohomologie halbeinfacher  
Matrizengruppen II, {\em Math. Z.} {\bf 92} (1966), 396--415.\\[-20pt]
\bibitem[Ill]{Ill} L.~Illusie, D\'eformations de groupes de  
Barsotti--Tate (d'apr\`es A. Grothendieck), in {\em S\'eminaire sur les  
pinceaux arithm\'etiques: la conjecture
de Mordell}, (Paris, 1983/84), Ast\'erisque {\bf 127} (1985), 151--198.\\[-20pt]
\bibitem[It\=o]{Ito} T. It\=o, Hasse invariants for some unitary
Shimura varieties, Algebraische Zahlentheorie (June 19th - June
25th, 2005), Mathematisches Forschungsinstitut Oberwolfach, Report
No. 28/2005, 1565--1568.\\[-20pt]
\bibitem[Ko]{Ko} R. Kottwitz, Points on some Shimura varieties over finite
fields, {\em J.~Amer.~Math.~Soc.}\ {\bf5} (1992), no. 2, 373--444.\\[-20pt]
\bibitem[Kr]{Kr} H.~Kraft, Kommutative algebraische p-Gruppen (mit  
Anwendungen auf
p-divisible Gruppen und abelsche Variet\"aten), manuscript 86 pages,  
Univ. Bonn, 1975.\\[-20pt]
\bibitem[LMB]{LMB} G. Laumon, L. Moret--Bailly, {\em Champs alg\'ebriques},
Ergebnisse der Mathematik und ihrer Grenzgebiete. 3. Folge, {\bf
39}, Springer-Verlag, Berlin, 2000.\\[-20pt]
\bibitem[Ma]{Ma} Y.~Manin, The theory of commutative formal
groups in finite caracteristic, {\em Russian Math. Surv.} {\bf 18}
(1963), no. 6, 1--83.\\[-20pt]
\bibitem[MFK]{MFK} D.~Mumford, J.~Fogarty, F.~Kirwan, {\em Geometric invariant
theory}, \emph{Third edition}, Ergebnisse der Mathematik und ihrer
Grenzgebiete. 2. Folge, {\bf 34}, Springer-Verlag, Berlin, 1994.\\[-20pt]
\bibitem[Mo]{Mo} B.~Moonen, Group schemes with additional structures  
and Weyl group
cosets, in {\em Moduli of Abelian Varieties} (Texel Island, 1999),  
255--298, Progr. Math., {\bf 195}, Birkh\"auser, Basel, 2001.\\[-20pt]
\bibitem[MW]{MW} B.~Moonen, T.~Wedhorn, Discrete invariants of  
varieties in positive
characteristic, {\em Int.~Math.~Res.~Not.} 2004, {\bf 72}, 3855--3903.\\[-20pt]
\bibitem[Oo1]{Oo1} F.~Oort, Newton polygons and formal groups:  
conjectures by Manin and Grothendieck,
{\em Ann. of Math.} (2)  {\bf 152}  (2000),  no. 1, 183--206.\\[-20pt]
\bibitem[Oo2]{Oo2} F.~Oort, A stratification of a moduli space of abelian
varieties, in {\em Moduli of Abelian Varieties} (Texel Island, 1999),  
345--416, Progr. Math., {\bf 195},
Birkh\"auser, Basel, 2001.\\[-20pt]
\bibitem[Tr1]{Tr1} C.~Traverso, $p$-divisible groups over fields,  
\emph{Symposia Mathematica}, Vol. {\bf XI} (Convegno di Algebra  
Commutativa, INDAM,
Rome, 1971), 45--65, Academic Press, London, 1973.\\[-20pt]
\bibitem[Tr2]{Tr2} C.~Traverso, Specializations of Barsotti--Tate  
groups, {\em Symposia Mathematica}, Vol. {\bf XXIV} (Sympos., INDAM,  
Rome, 1979),  1--21, Academic Press, London-New York, 1981.\\[-20pt]
\bibitem[Va1]{Va1} A.~Vasiu, Crystalline boundedness principle,
{\em Ann. Sci. \'Ecole Norm. Sup.} {\bf 39} (2006), no. 2,
245--300.\\[-20pt]
\bibitem[Va2]{Va2} A.~Vasiu, Level $m$ stratifications of versal  
deformations of
$p$-divisible groups, {\em J. Alg. Geom.} {\bf 17} (2008), no. 4, 599--641.\\[-20pt]
\bibitem[Va3]{Va3} A.~Vasiu, Mod $p$ classification of Shimura  
$F$-crystals, 63
pages, to appear in \emph{Math. Nachr.}
\\[-20pt]
\bibitem[Va4]{Va4} A.~Vasiu, Manin problems for Shimura varieties of  
Hodge type, 41
pages, available at
http://arxiv.org/abs/math/0209410.
\\[-20pt]
\bibitem[Wd]{Wd} T.~Wedhorn, The dimension of Oort strata of Shimura  
varieties of
PEL-type, in {\em Moduli of abelian varieties} (Texel Island, 1999),  
441--471, Progr. Math., {\bf 195},
Birkh\"auser, Basel, 2001.\\[-20pt]
\bibitem[Zi]{Zi} Th.~Zink, On the slope filtration, {\em Duke Math.  
J.} {\bf 109} (2001), no. 1, 79--95.

\end{thebibliography}
\end{document}